\documentclass[a4paper,11pt]{amsart}
\usepackage{amssymb}
\usepackage{amscd}
\def\Box{\square}

\def\frak{\mathfrak}
\def\Bbb{\mathbb}
\def\Cal{\mathcal}

\def\endrk{\hbox{$|\!\!|\!\!|\!\!|\!\!|\!\!|\!\!|$}}

\newcommand{\hook}{\raisebox{-0.35ex}{\makebox[0.6em][r]
{\scriptsize $-$}}\hspace{-0.15em}\raisebox{0.25ex}{\makebox[0.4em][l]{\tiny
 $|$}}}

\newcommand{\newc}{\newcommand}

\newcommand{\Sc}{\operatorname{Sc}}

\let\ccdot\cdot
\def\cdot{\hbox to 2.5pt{\hss$\ccdot$\hss}}

\newcommand{\om}{\omega}
\renewcommand{\phi}{\varphi}

\newcommand{\si}{\sigma}

\newcommand{\Om}{\Omega}

\newc{\aI}{\mbox{\boldmath{$ I$}}}
\newc{\aR}{\mbox{\boldmath{$ R$}}}
\newc{\aDeR}{\mbox{\boldmath{$ U$}}_B{}^P{}_C{}^Q}
\newc{\al}{\mbox{\boldmath$ \Delta$}}             
\newc{\nda}{\mbox{\boldmath$ \nabla$}}
\newc{\ad}{\mbox{\boldmath$ d$}}
\newc{\da}{\mbox{\boldmath$ \delta$}}
\newc{\aK}{\mbox{\boldmath{$ K$}}}
\newc{\aL}{\mbox{\boldmath{$ L$}}}

\newcommand{\g}{g^o}  

\newtheorem{theorem}{Theorem}[section]

\newtheorem{proposition}[theorem]{Proposition}
\newtheorem{corollary}[theorem]{Corollary}

\newcommand{\cN}{{\Cal N}}

\newcommand{\ce}{{\Cal E}}

\newcommand{\cq}{{\Cal Q}}

\newcommand{\cL}{{\Cal L}}
\newcommand{\cT}{{\Cal T}}

\newcommand{\ct}{{\Cal T}}

\newcommand{\bV}{{\Bbb V}}

\newcommand{\bS}{{\Bbb S}}

\newcommand{\nd}{\nabla}

\newcommand{\Up}{\Upsilon}

\newcommand{\Ric}{\operatorname{Ric}}

\newcommand{\nn}[1]{(\ref{#1})}

\newcommand{\bT}{\mathcal{T}} %{\mbox{{$\Bbb T$}}}

\def\W{\mathbb{W}}
\def\X{\mathbb{X}}
\def\Y{\mathbb{Y}}
\def\Z{\mathbb{Z}}

\newcommand{\bg}{\mbox{\boldmath{$ g$}}}
\newcommand{\V}{{\mbox{\sf P}}}                   
\newcommand{\J}{{\mbox{\sf J}}}

\newc{\strutdd}{\rule{0mm}{5mm}}

\newcommand{\lpl}                         
{\mbox{$
\begin{picture}(12.7,8)(-.5,-1)
\put(2,0.2){$+$}
\put(6.2,2.8){\oval(8,8)[l]}
\end{picture}$}}

\usepackage{ifthen}

\newc{\tensor}[1]{#1}

\newc{\Mvariable}[1]{\mbox{#1}}

\newc{\down}[1]{{}_{
\ifthenelse{\equal{#1}{;}}{|}{#1}}}

\newc{\up}[1]{{}^{#1}}
\newc{\C}{C}

\newc{\JulyStrut}{\rule{0mm}{6mm}}
\newc{\midtenPan}{\mbox{\sf S}}
\newc{\midten}{\mbox{\sf T}}
\newc{\midtenEi}{\mbox{\sf U}}
\newc{\ATen}{\mbox{\sf E}}
\newc{\BTen}{\mbox{\sf F}}
\newc{\CTen}{\mbox{\sf G}}
\newcommand{\RR}{\mathbb{R}}

\renewcommand{\S}{\Sigma}

\def\sideremark#1{\ifvmode\leavevmode\fi\vadjust{\vbox to0pt{\vss% the remark
 \hbox to 0pt{\hskip\hsize\hskip1em%                          will appear only
 \vbox{\hsize3cm\tiny\raggedright\pretolerance10000%          on the side
 \noindent #1\hfill}\hss}\vbox to8pt{\vfil}\vss}}}%
\begin{document}
%\begin{sloppypar}
\renewcommand{\today}{} \title{A class of compact Poincar\'e-Einstein
  manifolds:
 properties and construction}

\author{A. Rod Gover and Felipe Leitner}

\address{Department of Mathematics\\
  The University of Auckland\\
  Private Bag 92019\\
  Auckland 1\\
  New Zealand} \email{gover@math.auckland.ac.nz}

\address{ Universitat Stuttgart\\
Institut für Geometrie und Topologie\\
 Mathematics Department\\
 Pfaffenwaldring 57\\ 
D-70550 Stuttgart}\email{leitner@mathematik.uni-stuttgart.de}

\vspace{10pt}

\thanks{ARG gratefully acknowledges support from the Royal Society of
  New Zealand via Marsden Grant no.\ 06-UOA-029. L was supported by a
  research fellowship of the Deutsche Forschungsgemeinschaft which
  enabled him to visit the University of Auckland. FL is grateful to
  the Department of Mathematics at the UoA for hospitality. }

\renewcommand{\arraystretch}{1}
\maketitle
\renewcommand{\arraystretch}{1.5}
\pagestyle{myheadings}
\markboth{Gover \& Leitner}{Compact PE metrics}

\begin{abstract} 
  We develop a geometric and explicit construction principle that
  generates classes of Poincar\'e-Einstein manifolds, and more
  generally almost Einstein manifolds. Almost Einstein manifolds
  satisfy a generalisation of the Einstein condition; they are
  Einstein on an open dense subspace and, in general, have a conformal
  scale singularity set that is a conformal infinity for the Einstein
  metric.  In particular, the construction may be applied to yield
  families of compact Poincar\'e-Einstein manifolds, as well as
  classes of almost Einstein manifolds that are compact without
  boundary.  We obtain classification results which show that the
  construction essentially exhausts a class of almost Einstein (and
  Poincar\'e-Einstein) manifold. We develop the general theory 
of fixed conformal structures admitting
  multiple compatible almost Einstein structures. We also
  show that, in a class of cases, these are canonically related to a
  family of constant mean curvature totally umbillic embedded
  hypersurfaces.
\end{abstract}

\section{Introduction}

On a compact manifold $M^d$ with boundary, a metric $\g$ on the
interior $M^o$ of $M$ is said to be conformally compact if
$\g=s^{-2}g$ where $g$ is non-degenerate up to the boundary $\partial
M$ and $s$ is a defining function for $\partial M$. If also $\g$ is
Einstein then the structure is said to be Poincar\'e-Einstein (PE).
The term asymptotically hyperbolic Einstein is also used for these
structures which were introduced by Fefferman and Graham in
\cite{FGast} in connection with their ambient metric, and as a tool
for studying the conformal geometry of the boundary. The
relationship between that conformal structure and the geometry of the
Riemannian interior has recently been studied intensively using
spectral and scattering tools \cite{FGrQ,GrZ,Guill,Lspec,Wang}, and
related formal asymptotics \cite{Albin,CQY,FGrnew,GrJ}; some of these
treat structures which are suitably asymptotically PE. This
relationship is the geometric problem underlying the so called AdS/CFT
correspondence of String Theory \cite{Mal,Witten}.

The model of a PE manifold is the hyperbolic ball $\mathbb{H}^d$ and
this yields a large class of further examples via its quotients
$\Gamma\backslash \mathbb{H}^d$, where $\Gamma$ is a convex cocompact
discrete isometry group. There are also various specific examples known in
the literature, some of these are catalogued in the work \cite{AndAdS}
and among these are examples which are not locally conformally flat.
Anderson has rather general existence results in dimension 4
\cite{AndPresc}; this extends his earlier works \cite{AndL2,AndReg}. There
are now some rather powerful results available for constructing PE metrics
beginning with known ones; in particular through perturbing known
examples \cite{GrL,B,L}, and via Maskit combinations \cite{MP}.

Here we describe a product based construction principle for
constructing Poincar\'e-Einstein manifolds.  In fact
Poincar\'e-Einstein manifolds are a subclass of the so-called almost
Einstein (AE) manifolds as in \cite{GoPrague,GoIP,Go-al} and the
construction applies to this broader setting. Apart from the examples
of PE manifolds and Einstein manifolds, almost Einstein structures
arise naturally on certain classes of Fefferman spaces \cite{LMZ} (see
also \cite{CapGoF}) and in the constructions and classifications by
Derdzinski and Maschler of K\"ahler metrics which are ``almost
everywhere'' conformal to Einstein by a non-constant rescaling factor,
see e.g.  \cite{DML,DM} and references therein.  The construction we
develop enables the proliferation of a variety of classes of examples,
and in particular examples of compact PE manifolds and closed (i.e.\
compact without boundary) AE manifolds.  Warped products are well
established as a tool for the construction of Einstein metrics
\cite{Besse,ONeil}.  The construction here may be viewed as
generalising, or extending, a class of generalised warped products.
This extension is made transparent through the fundamental
relationship between (almost) Einstein structures and the prolonged
differential system captured by the conformal tractor connection.
Among the simplest examples, we obtain almost Einstein structures on
the sphere products $S^2\times S^1$ and $S^3\times S^1$, manifolds
well known to not admit Einstein metrics \cite{Besse}.

In \cite{GLsub} for appropriate products of Einstein manifolds we
describe the explicit construction of Poincar\'e-Einstein collar
metrics which have the given product as conformal infinity. One may
ask which of these may be compactified, that is, which may be embedded
into compact PE manifolds. This is essentially answered here in the
classification results of Section \ref{class}. We show in Corollary
\ref{CorBB} that a class of these collars may not be compactified. The
remainder, at least given some assumptions concerning simple
connectivity, will be seen to arise as conformal infinity
neighbourhoods of the examples of Section \ref{SM}, see Theorem
\ref{glob}.

All constructions in the following will be on Riemmanian manifolds
$(M,g)$ and Riemannian signature conformal manifolds $(M,c)$ (where
$c$ is an equivalence class of Riemannian metrics, each pair related
by rescaling by a positive function). Many of the statements hold
true in other signatures but to simplify the exposition we shall
simply avoid this point.  Recall that on $(M^d,g)$, with $d\geq 3$, the
Schouten tensor $P$ is the trace modification of the Ricci tensor
$\Ric$ which satisfies $\Ric =(d-2)P+gJ$, where $J$ is the metric
trace of $P$.  (Until further notice we shall restrict to manifolds of
dimension $d\geq 3$.)  We will say that $(M,g,s)$ (or really this modulo an equivalence relation, see Section \ref{AES}) is a {\em directed
almost Einstein}  structure if $s\in C^\infty (M)$ is a
non-trivial solution to the equation
\begin{equation}\label{prim}
trace-free(\nabla^g\nabla^g
s + sP^g) =0 .
\end{equation}

In \cite{GoIP} we see that Poincar\'e-Einstein manifolds are a special
class of directed AE manifolds.  On the other hand, from the following
result we see that a large class of directed AE manifolds are
necessarily boundary gluings of PE manifolds.
\begin{theorem} \cite{Go-al} \label{maina} Let $(M^d,g,s)$ be a directed almost
  Einstein structure with $g$ positive definite and $M$ connected.
  Then writing $\S$ for the zero set of $s$, $M\setminus \Sigma$ is an
  open dense set and on this $\g=s^{-2}g$ is Einstein with $\Ric=
  (d-1)Sg$. According
  to the sign of the constant $S=S(g,s)$, there are three cases: \\
  $\bullet$ If
  $S>0$ then the scale singularity set $\S$ is empty. \\
  $\bullet$ If $S=0$ then $\S$ is either empty or otherwise consists
  of isolated
  points and these points are critical points of the function $s$.\\
  $\bullet$ If $S<0$ then $\S$ is either empty or else is a totally
  umbillic smooth hypersurface.  In particular on closed (i.e.\ $M$ compact
  without boundary) $S<0$ almost Einstein manifolds, with $\S\neq
  \emptyset$, a constant rescaling of $s$ normalises $S$ to $-1$, and
  then $(M\setminus M^-)$ is a finite union of connected
  Poincar\'e-Einstein manifolds. Similar for $(M\setminus M^+)$. Here $M^\pm$ is the open submanifold where $s$ is positive, respectively, negative. 
\end{theorem}
\noindent Much more can be said about the scale singularity set $\S$
(the conformal infinity of the structure) and its relation to the
ambient geometry. This is one of the main directions of \cite{Go-al}.  Here
the central point is that, from the last item of the above,
constructing scalar negative directed AE manifolds is an effective
route to PE structures.

If $s$ solves \nn{prim} then so does $-s$, and where $s$ is
non-vanishing these solutions determine the same Einstein metric. We
shall say that a manifold $(M,g)$ is {\em almost Einstein} (AE) if it
admits a covering such that on each open set $U$ of the cover we have
that $(U,g,s_U)$ is directed almost Einstein and on overlaps $U\cap V$
we have either $s_U=s_V$ or $s_U=-s_V$. An AE structure will be said
to be {\em directable} if there is a directed AE structure that agrees
locally and up to sign with the given AE structure. 
Just as the standard conformal sphere yields a
model for all categories of directed almost Einstein structure
\cite{Go-al} we observe in Section \ref{models} that the corresponding
conformal structure on $\mathbb{RP}^d$ yields a model for AE
structures which are not directable.  These yield further examples
through the construction principle of Section \ref{wprod}, see the examples of Section
\ref{nond}.

The tractor connection, due to Thomas \cite{T}, is a natural and
conformally invariant vector bundle connection equivalent (see
\cite{CapGoTAMS}) to the normal conformal Cartan connection.  As we
shall review briefly in the next section, parallel sections of this
are in 1-1 correspondence with solutions of \nn{prim}; in fact more
than this, from the construction of the tractor connection in
\cite{BEG} one sees that its equation of parallel transport just {\em
  is} a prolonged system for \nn{prim} (cf.\ \cite{BEG}). {\em So a directed
  AE structure is a conformal manifold equipped with a parallel
  section of the (standard) tractor bundle.} More generally an AE
manifold is one which admits a parallel tractor up to a sign in the
obvious way. Using this, it follows that a simply connected AE
manifold necessarily admits a directed almost Einstein structure.

Each of the closed almost Einstein spaces we construct in Section \ref{wprod}
admits, locally at least, more than one linearly independent
parallel tractor; i.e., the conformal class of each admits at least two
non-trivially distinct  almost Einstein structures.
Following some preliminaries, in Sections \ref{crE} and \ref{crae} we
investigate the consequences of a fixed conformal structure admitting
two linearly independent parallel tractors. This is a natural
extension of the problem of conformally related Einstein structures.
As in that case, we see here that for closed manifolds the condition
of multiple almost Einstein structures is very restrictive.  
Let us write $\mathcal{K}_c$ for the vector space of parallel tractors
on a given conformal manifold $(M,c)$. The main result of Section
\ref{crE} is the following.
\begin{theorem} \label{nonitsure} Let $(M^d,c)$ be a closed
 Riemannian conformal space of dimension $d\geq 3$ with $\dim
 (\mathcal{K}_c )\geq 2$.  Then either
\begin{enumerate}
\item $(M,c)$ is the standard sphere
$(S^d,[g_{\rm rd}])$, or

\item for any $K\in\mathcal{K}_c\smallsetminus\{0\}$, it is necessarily
  the case that $S_K <0$ and $\Sigma_K$ is non-empty (and hence is a
  totally umbillic hypersurface in $(M,c)$).

\end{enumerate}
\end{theorem}
\noindent In the Theorem $S_K$ is the constant $S$, as in Theorem
\ref{maina}, for the almost Einstein structure $(M,c,s_K)$ where $s_K$
is determined by the parallel tractor $K$.  Note that, in particular,
excluding the round sphere, the conformal structure of a closed
Einstein manifold never admits a {\em proper} directed almost Einstein
structure, meaning one with a non-trivial conformal infinity.
Equivalently, and again excepting the sphere, on a closed manifold a
proper directed almost Einstein structure does not have an Einstein
metric in the conformal class. This is part of a picture which
indicates that proper almost Einstein structures have in a sense ``the
same status'' as Einstein structures.

It turns out that more can be said about the structure of the scale
singularity sets when $\dim (\mathcal{K}_c )\geq2$.
\begin{theorem}\label{min}
Let $(M^d,c)$ be a Riemannian conformal space of dimension $d\geq 3$.
If $K_1,K_2\in \mathcal{K}_c$ are linearly independent with
$S_{K_2}<0$ and $\S_{K_2}\neq \emptyset$, then on $M\setminus
\Sigma_{K_1}$ $\Sigma_{K_2}$ (is totally umbillic and) has constant
mean curvature with respect to the Einstein metric $\g$ determined by
$K_1$.
If the tractors $K_1$ and $K_2$ are orthogonal then $\Sigma_{K_2}$ is
minimal (and hence totally geodesic) with respect to the Einstein
metric $\g$.
\end{theorem}
\noindent We shall prove this in Section \ref{crae}.
Related to this we shall see that in fact there is a
generalised notion of mean curvature which extends to the conformal
infinity $\Sigma_{K_1}$ (in the case this is not empty) and this is
closely linked to the scalar curvature quantity $S_{K_1}$.

Given two Einstein manifolds $(M_1,g_1)$ and $(M_2,g_2)$, if there is
a suitable relationship between the scalar curvatures of each, then we
term the product $(M,g)=(M_1\times M_2,g_1\times g_2)$ a special
Einstein product. A key observation is that, for such products, the
conformal tractor connection is a simple restriction of the direct sum
connection. It is convenient for our purposes here to see this
directly via the definition and characterisation of the tractor
bundle, and this is treated in Proposition \ref{prodT}. (See also
\cite{LeitNCK} and \cite{Armcomplete} where  this is observed and
is seen to arise naturally in connection with decomposable conformal
holonomy. It should be noted that this result also follows
easily from the explicit Fefferman-Graham (ambient) metric
construction for special Einstein products from Theorem 2.1 of
\cite{GLsub} via the Fefferman-Graham metric/tractor relationship
established in \cite{CapGoamb}.) Using that the tractor connection is
essentially a direct sum it follows that a special Einstein product is
almost Einstein if and only one of the factors (say $(M_2,g_2)$)
admits an almost Einstein structure linearly independent of the
declared Einstein structure ($g_2$).  Here the ``linear independence''
refers to the corresponding parallel tractors. This result, Theorem
\ref{prodAE}, is the basis of the construction principle. Also
important in this respect is that if one of the factors has a proper
almost Einstein structure then the corresponding AE structure on the
product is also proper and the singularity set is easily understood,
see Corollary \ref{singprod}. Away from the singularity set the
Einstein metric $\g$ (notation as in Theorem \ref{maina}) is a
generalised warped product metric of the form \nn{dcneg}. 
Thus
the construction here of almost Einstein manifolds may be viewed as an
extension or generalisation of the construction of Einstein manifolds
by (suitable) warped products of this form.

The construction via Theorem \ref{prodAE} of almost Einstein manifolds
that are proper and closed is particularly interesting since then we
also obtain compact Poincar\'e-Einstein spaces via Theorem
\ref{maina}.  On the other hand the possibilities are limited by
Theorem \ref{nonitsure}; for the construction of directed AE special
Einstein products one of the factors must be a sphere.  Thus the
classes of examples developed in Sections\ \ref{SM} through to Section
\ref{exs1} are exhaustive.  Section \ref{wprod} concludes with some
non-directed almost Einstein structures; note that example \ref{nond}
with $m_1=1$ shows that any compact Einstein manifold arises as the
conformal infinity of a closed almost Einstein manifold. The examples
\ref{SM} give closed conformal manifolds admitting multiple compatible
directed AE structures. Generically for these examples the scale
singularity sets of the linearly independent parallel tractors
intersect. When this happens, if we pick one of these to determine the
Einstein metric $\g$ then any other gives a totally umbillic
hypersurface (of constant mean curvature in the sense Theorem
\ref{min}) that meets the conformal infinity of $\g$ in a totally
umbillic hypersurface (or in isolated double points), see Proposition
\ref{minmore}. Since this is a way of generating minimal submanifolds
that meet conformal infinity it seems that these may be interesting
models for branes and the study of their asymptotics (cf.\ \cite{AM,GW}).

\section{Conformal Geometry, tractor calculus and almost Einstein structures} 

\subsection{The conformal tractor bundle}\label{trba}
Let $M^d$ be a smooth manifold. Until further notice the dimension $d$
is at least 3.  It will be convenient to use some standard structures
from conformal geometry, further details and background may be found
in \cite{CapGoamb,GoPetCMP}.  Recall that a (Riemannian) {\em
conformal structure\/} on $M$ is a smooth ray subbundle $\cq\subset
S^2T^*M$ whose fibre over $p\in M$ consists of conformally related
positive definite metrics at the point $p$. Sections of $\cq$ are
metrics $g$ on $M$. So we may equivalently view the conformal
structure as the equivalence class $c=[g]$ of these conformally
related metrics.  The principal bundle $\pi:\cq\to M$ has structure
group $\Bbb R_+$, and so each representation $\mathbb{R}_+ \ni
t\mapsto t^{-w/2}\in {\rm End}(\mathbb{R})$ induces a natural line
bundle on $ (M,c)$ that we term the conformal density bundle $\ce[w]$.
We write $\ce^a$ as an abstract index notation for the tangent bundle
$TM$ and similarly $\ce_a$ is an alternative notation for the
cotangent bundle $T^*M$.  In general each vector bundle and its space
of smooth sections will be denoted by the same notation.

We write $\bg$ for the {\em conformal metric}, that is the
tautological section of $S^2T^*M[2]:= S^2T^*M\otimes \ce[2]$
determined by the conformal structure. This will be henceforth used to
identify $TM$ with $T^*M[2]$.  For example, with these conventions the
Laplacian $ \Delta$ is given by $\Delta=-\bg^{ab}\nd_a\nd_b=
-\nd^b\nd_b\,$ where $\nabla$ (or sometimes we will write $\nabla^g$)
is the Levi-Civita connection for some choice of metric $g$ from the
conformal class.  Note $\ce[w]$ is trivialised by a choice of metric,
and we write $\nd$ (or again sometimes $\nabla^g$) for the connection
corresponding to this trivialisation.  It follows immediately that the
(coupled) connection $ \nd_a$ preserves the conformal metric and so
the canonical conformal volume form (of conformal weight $n$) is also
parallel for $\nd_a$.  The conformal metric $\bg$ and its inverse will
henceforth be the default object used to contract indices on tensors
even when we have fixed a metric from the conformal class.

The Riemann curvature tensor $R_{ab}{}^{c}{}_d$ is given by
$$(\nabla_a\nabla_b-\nabla_b\nabla_a)V^c=R_{ab}{}^{c}{}_d V^d 
\quad\text{
where} \quad \ V^c\in \ce^c.$$
This can be decomposed into the totally trace-free {\em Weyl curvature}
$C_{abcd}$ and the symmetric
Schouten tensor $P_{ab}$ according to
$$
R_{abcd}=C_{abcd}+2g_{c[a}P_{b]d}+2g_{d[b}P_{a]c}.
$$
Changing to a conformally related metric $\widehat{g}=e^{2\om}g$
($\om$ a smooth function) one finds the Weyl tensor is conformally
invariant $\widehat{C}_{ab}{}^c{}_d=C_{ab}{}^c{}_d$, whereas the
Schouten tensor transforms to
\begin{equation}\label{Rhotrans}
\textstyle \widehat{\V}_{ab}=\V_{ab}-\nd_a \Up_b +\Up_a\Up_b
-\frac{1}{2} \Up^c\Up_c \bg_{ab} ,
\end{equation}
where $\Upsilon:=d\om$.

For a given choice of metric $g$, the (standard) {\em tractor} bundle
$\ct$, or using an abstract index notation $\ce^A$, may be identified
with a direct sum
$$
\ct=\ce^A \stackrel{g}{=} \ce[1]\oplus\ce_a[1]\oplus\ce[-1] . 
$$ Thus a section $U$ of $\ct$ may be identified with a triple
$(\si,~\mu_a,~\rho)$; when it is understood that a metric $g$ has been
fixed we will write simply $U^A= (\si,~\mu_a,~\rho)$.  Conformally
transforming, as above, gives a different isomorphism, which is
related to the previous by the transformation formula
\begin{equation}\label{transf}
\widehat{(\si,\mu_b,\rho)}=(\si,\mu_b+\si\Up_b,\rho-\bg^{bc}\Up_b\mu_c-
\tfrac{1}{2}\si\bg^{bc}\Up_b\Up_c).  
\end{equation}
It follows that $\bT$ invariantly has a composition structure
$$
\bT=\ce[1]\lpl \ce_a[1]\lpl \ce[-1];$$ 
$\ce[-1]$ may be identified
with a subbundle of $\bT$ and $\ce_a[1]$ is a subbundle of the
quotient bundle $\bT/\ce[-1]$.  From \nn{transf} one also sees that
the map $\ce^A \to \ce[1]$ is conformally invariant and may be
regarded as a preferred section $X_A$ of $\ce_{A}[1]$ so that, with
$U^A$ as above, we have $\sigma = U^{A}X_{A}$.  This section also
describes the invariant injection $\ce[-1] \to \ce_{A}$ according to
$\rho \mapsto \rho X_A$. In computations, it is often useful to
introduce the remaining `projectors' $Z_{Aa}:\ce^A\to \ce_a[1]$ and
$Y_A:\ce^A\to \ce[-1]$.  which are determined by a choice of scale.
For calculations it is useful to view these as sections
$Z_{Aa}\in\ce_{Aa}[1]$ and $Y_A\in\ce_A[-1]$, where
$\ce_{Aa}[w]=\ce_A\otimes\ce_a\otimes\ce[w]$, and so forth.

We describe any tensor product (or symmetric tensor product etcetera)
of the tractor bundle and its dual as tractor bundles.  If such a
bundle is tensored with some bundle of densities $\ce[w]$ then we
shall describe the result as a {\em weighted tractor bundle}. 
Repeated (abstract) tractor indices indicate a contraction, just as for
tensor indices.

The bundle $\ce^A$ carries an invariant signature $(d+1,1)$ {\em
  tractor metric} $h_{AB}$, and a connection $\nd_a$ which
  preserves this.  For $ U^A$ as above, this metric is determined by 
\begin{equation}\label{trmet}
h_{AB}U^AU^B =2\rho \sigma +\mu^a\mu_a.
\end{equation}
As a point on notation, we may also write $h(U,U)$ for the expression
in the display.  The tractor metric will be used to raise and lower
indices without further mention.  In terms of the metric $g$ from the
conformal class, the tractor connection is given explicitly by
the following formula for $\nd_a U^B$:
\begin{equation}\label{ndef}
    \nd_a \left(\begin{array}{c} \sigma \\
  \mu^b \\ \rho \end{array}\right)
  = \left(\begin{array}{c} \nd_a \sigma - \mu_a \\
  \nd_a\mu^b + \delta_a{}^b \rho +
  P_{a}{}^b \sigma \\ \nd_a\rho - P_{ac} \mu^c
  \end{array} \right) .
\end{equation}
This gives the normal tractor connection, in all statements below the
tractor connection may be assumed to be this one, unless otherwise
stated.  Of course this may be extended to a connection on any tractor
bundle in the obvious way. The use of the same symbol $\nd$ as for the
Levi-Civita connection is intentional. More generally, we shall use
$\nd$ to mean the coupled Levi-Civita-tractor connection: this enables
us, for example, to apply $\nd $ to weighted tractor bundles or
tensor-tractor. Although in this case it is not conformally invariant
it enables us to, for example, compute the covariant derivative of the
tractor projectors $X$, $Y$ and $Z$.

As discussed in \cite{BEG},  there is  an invariant second order
operator between weighted tractor bundles due to T.Y. Thomas,
$$
D_A:\ct^\ast[w]\to\ct^\ast[w-1] , \quad \quad \ct^*~\mbox{ any tractor bundle},
$$
given by 
\begin{equation}\label{Dform}
D_A V:=(d+2w-2)w Y_A V+ (d+2w-2)Z_{}^{a}\nabla_a V + X^A(\Delta-w J)V.
\end{equation} 
For an invariant construction of this see \cite{GoSrni99,CapGoluminy}.

\subsection{Almost Einstein structures} \label{AES}
As above we work on Riemannian manifolds $(M^d,g)$ with the dimension
$d\geq 3$. The Schouten tensor $P$ (or $P^{g}$), introduced
earlier, is related to the Ricci tensor by $ \Ric=(d-2)P+J g $, where
recall $J$ is the conformal metric trace of $P^g$. Calculating the conformal
transformation of $P$, it follows easily that the metric $g$ is
conformally Einstein if and only if there is a nowhere zero solution
$s\in C^\infty(M)$ to the equation \nn{prim} $trace-free(\nabla\nabla
s + P s)=0$.

Note that the equation \nn{prim} is conformally covariant. 
 If we write
$A(g,s):=trace-free(\nabla\nabla s + P s) $ then $e^\om A(g,s)=
A(e^{2\om}g,e^{\om}s)$, $\om \in C^\infty (M)$. By a directed almost
Einstein structure we mean $(M,g,s)$, where $s$ is a solution of
\nn{prim}, but modulo the equivalence relation $(g,s)\sim
(e^{2\om}g,e^{\om}s)$; more precisely it is an equivalence class
$[(M,g,s)]$ where the equivalence relation is $(M,g,s)\sim
(M,e^{2\om}g,e^{\om}s)$ and $s$ is a solution of \nn{prim}.

The conformal covariance of the equation \nn{prim} is equivalent to
conformal invariance, if we replace the function $s$ by a conformal
density $\si\in \ce[1]$. The equation 
\begin{equation}\label{primc}
trace-free(\nd_a\nd_b \si+P_{ab} \si)=0
\end{equation} 
is conformally invariant. 
This suggests a strong link to conformal
geometry.  As observed in \cite{GoPrague} it follows easily from
the definition \nn{ndef} of the tractor connection above that we have
the following.
\begin{theorem}\label{key}
A directed almost Einstein structure is a conformal manifold $(M,c)$
 equipped with a parallel (standard) tractor $I\neq 0$. The mapping
 from non-trivial solutions of \nn{primc} to parallel tractors is by
 $\si\mapsto \frac{1}{d}D\si$ with inverse $I\mapsto \si:=h(I,X)$, and $\si$ is
 non-vanishing on an open dense set $M\setminus \Sigma$. On this set
 $\g:=\si^{-2}\bg$ is Einstein.
\end{theorem}
The point here, as touched on earlier, is the parallel transport
equation for the tractor connection \nn{ndef} is simply a prolonged
system derived from (and essentially equivalent to) \nn{primc}
\cite{BEG,Go-al}.  So a directed almost Einstein manifold is exactly a
conformal manifold $(M,c)$ equipped with a parallel standard tractor $I$ and
we write $(M,c,I)$ to indicate this.
The set $\Sigma$, where the almost Einstein ``scale'' $\si:=h(X,I)$
vanishes, shall be called the {\em scale singularity set}.  

On an almost Einstein manifold
$(M,c,I)$ (or $(M,g,s)$) we shall write $S_I$ or $S(\si)$ (or even $
S(g,s)$) as a shorthand for $-h(I,I)$. This may be viewed as a scalar
curvature quantity for the structure, since off $\Sigma$ we have
$S(\si)= \frac{\Sc^{\g}}{d(d-1)}$, and this is the quantity $S$ in
Theorem \ref{maina}.  In this language that theorem is as follows.
\begin{theorem}\label{classthm} \cite{Go-al} 
  Let $(M,c,I)$ be a Riemannian signature almost Einstein structure
  and $\si=h(X,I)$.  If $S(\si) >0 $ then $\Sigma$ is empty and
  $(M,\si^{-2}\bg)$ is Einstein with positive scalar curvature; If
  $S(\si)=0$ then $\Sigma$ is either empty or consists of isolated
  points $p$ where $j^1_p\si=0$, and $(M\setminus \Sigma,\si^{-2}\bg)$
  is Ricci-flat; if $S(\si) < 0$ then the scale singularity set
  $\Sigma$ is either empty or else is a totally umbillic smooth hypersurface,
  and $(M\setminus \Sigma, \si^{-2} \bg)$ is Einstein of negative
  scalar curvature.
\end{theorem}
\noindent
The Theorem here is a straightforward consequence of the formulae
\nn{ndef} and \nn{trmet}. There is a detailed treatment of this and
surrounding results in \cite{Go-al}.  Note that if $\S\neq \emptyset$
then $\S$ is a conformal infinity for the ambient Einstein metric $\g$
in the usual sense; there is a metric $g$ on the manifold $M$ so that
$\g=s^{-2}g$. In this situation it is well known that if $M$ is
complete then $(M\setminus \S,\g)$ is geodesically complete. In the
case that $M$ is compact it is traditional to describe $(M\setminus
\S,\g)$ as conformally compact.  If also $S(\si)<0$ then the smooth
hypersurface $\S$ inherits the conformal structure $(\S,[g|_\S])$.

 In terms of the tractor picture to say that we have an almost
Einstein manifold that is not necessarily directed means that on each
open set $U$ of some cover we have $(M,c,I_U)$, with $I_U$ parallel,
and on overlaps $U\cap V$ we have $I_U=\pm I_V$.  We will write
$(M,c,[I])$, or simply $[I]$ if the conformal structure is fixed, to
denote such an AE structure.

\subsection{The Models}\label{models}

On the $d$-sphere, with its standard conformal structure, the standard
tractor bundle may be trivialised globally by parallel sections. Thus
we have the following.
\begin{proposition}\label{modp}
  The $d$-sphere, with its standard conformal structure, admits a
  $(d+2)$-dimensional space of compatible directed almost Einstein
  structures.
\end{proposition}
In fact for each $S\in \mathbb{R}$ there is an almost Einstein
structure $I$ on $\bS^d$ with $S_I=S$, this is Proposition 5.1 in
\cite{Go-al} and the case of the sphere as a model for directed almost
Einstein structures is treated in detail there.

The standard conformal structure on the sphere is recovered as
follows.  Consider a $(d+2)$-dimensional real vector space $\bV$
equipped with a non-degenerate bilinear form $\mathcal{H}$ of
signature $(d+1,1)$. The {\em null cone} $\cN$ of zero-length vectors
form a quadratic variety.  Let us write $\cN_+$ for the forward part
of $\cN\setminus \{0 \}$.  Under the ray projectivisation of $\bV$ the
forward cone $\cN_+$ is mapped to a quadric in $\mathbb{P}_+(\bV) \cong
\bS^{d+1}$.  This image is topologically a sphere $\bS^d$ and we will
write $\pi$ for the submersion $\cN_+\to \bS^d$. Each point $p\in
\cN_+$ determines a positive definite inner product on
$T_{x=\pi{p}}\bS^d$ by $g_x(u,v)=\mathcal{H}_p(u',v')$ where $u',v'\in
T_p\cN_+$ are lifts of $u,v\in T_x\bS^d$. For a given vector $u\in T_x
\bS^{d}$ two lifts to $p\in \cN_+$ differ by a vertical vector field.
Since any vertical vector is orthogonal (with respect to $\mathcal{H}$) to
the cone it follows that $g_x$ is independent of the choices of lifts.
Clearly then, each section of $\pi$ determines a metric on $\bS$ and
by construction this is smooth if the section is. (Evidently the
metric agrees with the pull-back of $\mathcal{H}$ via the section
concerned.) Now, viewed as a metric on $T\mathbb{R}^{d+2}$,
$\mathcal{H}$ is homogeneous of degree 2 with respect to the standard
Euler vector field $E$ on $\bV$. That is $\cL_E \mathcal{H}=2
\mathcal{H}$, where $\cL$ denotes the Lie derivative. In particular,
this holds on the cone, and we note also that the cone is generated by $E$.

As explained in \cite{Go-al} (and cf.\ \cite{CapGoamb}), the tractor
bundle on $\bS^d$ is then naturally identified with $T\bV|_{\cN_+}
/\sim$. Here the equivalence relation ``$\sim$'' is as follows: for
$p,q\in \cN_+$ we have $T_p\sim T_q$ if and only if $p$ and $q$ lie in
the same null ray through the origin and $T_p$ is parallel to $T_q$
according to the usual parallelism of $\bV$ as an affine space. The
tractor connection is induced from this same parallelism, and the
usual tractor metric $h$ is induced from $\mathcal{H}$.

It follows that each vector $I$ in $\bV$ determines a  parallel
tractor field on $\bS^d$. In particular if we fix a choice $I$ with
$|I|^2=-1$ the we obtain a standard round metric on the sphere; we
shall write $(S^d,g)$ for this. This can be viewed as arising from the
section of $\cN_+$ given by the intersection with the hyperplane
$I_AX^A=1$, where $X^A$ are standard coordinates on $\bV\cong
\mathbb{R}^{d+2}$.  Let us write $I^\perp$ for the vectors in $\bV$
orthogonal to $I$ according to $\mathcal{H}$ and then also use the
same notation for the corresponding tractors. It is easily verified
that the tractors $K\in I^\perp$ correspond to solutions of \nn{primc}
that are odd with respect to $(S^d,g)$; $\alpha:=h(K,X)$ takes the
opposite sign at antipodal points. Let us now form from $(S^d,g)$ the
{\em unit $\mathbb{RP}^d$} which we shall denote by
$(\mathbb{RP}^d,g)$ by taking an antipodal identification of
$(S^d,g)$. It is clear that each $K\in I^\perp $ descends to an almost
Einstein structure $[K]$ on $(\mathbb{RP}^n,g)$ that is not directed.
Thus we have the following.
\begin{proposition}\label{modpnd}$\mathbb{RP}^d$
  with its standard conformal structure, admits an Einstein metric and
  also a $(d+1)$-dimensional space of compatible almost Einstein
  structures that are not directed.
\end{proposition}

\subsection{Conformal Killing vectors and tractors}

\newcommand{\K}{\mathbb{K}}

Via Theorem \ref{key} the equation of tractor parallel transport  
captures the  basic equation \nn{prim} in a geometric form that  
plays a central role in our discussions below. 
The  prolonged system corresponding to conformal
Killing vectors is also handled well by the tractor machinery. 
First note that (as observed in \cite{Esrni}) 
the operator $M_{A}{}^a:\ce_a[w]\to \ce_A[w-1]$ given, in a conformal
scale by
$$
u_a\mapsto (d+w-2)Z_A{}^a u_a -X_A\nd^a u_a,
$$
is conformally invariant.  
Then by Lemma 2.1 in \cite{powerslap} we have that solutions $k^a$ of
the conformal Killing equation
$$
\nd_{(a}k_{b)_0}=0
$$
are in 1-1 correspondence with fields $K_B\in\ce_B[1]$ 
such that  
\begin{equation}\label{cksttr}
D_{A}K_{B}= \K_{AB}
\end{equation} 
where $\K_{AB}$ is skew (i.e. it is a {\em 2-form tractor}). 
The correspondence is given by 
$$ k^a\mapsto \frac{1}{d}M_{Ba}k^a \quad \mbox{with inverse} \quad K_{B}
\mapsto Z^{Ba}K_{B}.
$$  

Using the above (or otherwise \cite{Capinf,GoSilKill}), the conformal
Killing condition is also captured by an equation on the skew part
$\K_{AB}$ of $D_A K_B$ as follows.
\begin{proposition}\cite[Proposition 3.1]{powerslap}\label{cKilltractor}
Solutions $k^a$ of the conformal Killing equation
$$
\nd_{(a}k_{b)_0}=0
$$
are in 1-1 correspondence with solutions $\K_{DE}\in \ce_{[DE]}$ of the 
equation
\begin{equation}\label{cktreq}
\nd_a \K_{DE}-\K_{AB}X^AZ^{Ba}\Omega_{abDE}=0.
\end{equation}
The correspondence is given by 
$$ k^a\mapsto \frac{1}{d^2} D_{[A}M_{B]a}k^a \quad \mbox{with inverse} 
\quad \K_{AB}
\mapsto X^AZ^{Ba} \K_{AB}.
$$  

If $\K_{AB}$ is a parallel adjoint tractor then $k^a:=X^AZ^{Ba}
\K_{AB}$ is a conformal Killing vector field and hence, from \nn{cktreq}, satisfies $k^a\Omega_{abDE}=0 $.
\end{proposition}

Given $k_a\in \ce[2]$ the 2-form tractor $\K_{AB}$ is given by
\begin{equation}\label{KKform}
\K_{AB}\stackrel{g}{=}\left( \begin{array}{ccc} & k_a &\\
\nabla_{[a}k_{b]} && \frac{1}{d}\nabla^ck_c \\
&\frac{1}{d}( \nabla^c\nabla_{(c} k_{a)_0} -\nabla_{a}\nabla^c k_c ) -P_a{}^ck_c& \end{array}
\right)
\end{equation}
where from the top-left to right-bottom the entries are the
coefficients of $\Y_{AB}=Y_{[A}Z_{B]}{}^a$, $\Z_{AB}^{ab}=
Z_{[A}{}^aZ_{B]}{}^b$, $\W_{AB}=X_{[A}Y_{B]}$,
$\X_{AB}=X_{[A}Z_{B]}{}^a$.

We say a conformal vector field on a space $(M,c)$ is {\em essential}
if it is not a Killing vector with respect to any metric $\tilde{g}\in c$ 
on $M$, i.e.
\[ L_V\tilde{g}=\frac{2}{d}div_{\tilde{g}}(V) \tilde{g}\neq 0\qquad
\mbox{for\ any}\ \tilde{g}\in c\ .\] Note that every essential
conformal vector field $V$ necessarily has a zero at some point; a
conformal Killing vector $V$ without a zero is Killing for the
conformally changed metric $\frac{1}{g(V,V)}g$.

We call a
conformal Killing vector field $k$ on $(M,g)$ a {\em conformal gradient} if
there exists a smooth function $\psi$ on $M$ such that
$k=grad_g(\psi)$.  Note that this is not a conformally invariant
property; $k$ is not a gradient for other metrics in the conformal
class.
\begin{theorem}\label{kzeros} On a 
  Riemannian manifold, of dimension at least 3, a conformal gradient
  with at least one zero is an essential conformal vector field.
 \end{theorem}
\noindent{\bf Proof:} If $k$ is a gradient then in any scale we have
$\nabla_{[a}k_{b]}=u_{[a}k_{b]}$ for some 1-form $u$. If we assume
also that $k$ is not essential and for $p\in M$ we have $k(p)=0$ then
it follows that the tractor 2-form $\mathbb{K}_{AB}$ (of \nn{KKform})
is zero at $p$. But the left hand side of \nn{cktreq} determines a
connection (which is a modification of the normal tractor
connection) for which this tractor $\mathbb{K}$ is parallel. So
$\mathbb{K}$ is zero everywhere and hence so is $k$.  \quad $\Box$\\
In fact this result is well known and also holds in dimension 2
\cite{KRI} (although we shall use it only in higher dimensions);
we include it here for completeness since, as seen, it follows easily
from Proposition \ref{cKilltractor}.

\subsection{Multiple almost Einstein structures}\label{crE}

Here we want to consider the situation where a given {\em fixed
  conformal structure} admits linearly independent parallel standard
tractors, or equivalently, linearly independent solutions of
\nn{primc}. Note that in general these are not ``conformally
related'', as two solutions will in general have different zero sets.
Of course if $\si_1$ and $\si_2$ are two such solutions, then away
from the union of the two zero sets $\si_1^{-2}\bg$ and
$\si_2^{-2}\bg$ are conformally related Einstein metrics.

It is well known that a pair of conformally related Einstein metrics
determines a conformal gradient field \cite{Brinkmann}. We observe here that
this generalises (adapting and extending  an
argument from Theorem 2.4 of \cite{powerslap}).
\begin{theorem}\label{2eincgrad}
  If $g_1 = (\si_1)^{-2}\bg$ is an Einstein metric, then
$(M,[g_1],\si_2)$ is an almost Einstein structure if and only if the
vector field
\begin{equation}\label{kdef}
 k^a:=\si_1\nd^a\si_2 -\si_2\nd^a \si_1
\end{equation}
is a conformal gradient field with respect to $g_1$. 
\end{theorem}
\noindent{\bf Proof:} 
 $\Rightarrow:$ Since $\si_1$ and $\si_2$ are almost
Einstein scales it follows from Theorem \ref{key} that
$$
I_1^A := \frac{1}{n} D^A \si_1 \quad \mbox{ and } I_2^A:= \frac{1}{n}D^A\si_2 
$$
are both parallel for the tractor connection. Thus 
\begin{equation}\label{kpar}
  \K^{AB}:= I_1^AI_2^B- I_1^BI_2^A \stackrel{g}{=}
  \left( \begin{array}{c}  \si_1\nabla^a\si_2- \si_2\nabla^a\si_1 \\
      (\nabla_{[a}\si_1)\nabla_{b]}\si_2 \phantom{wordortwo} \frac{1}{d} (\si_2\Delta \si_1- \si_1\Delta \si_2) \\
      \frac{1}{d}\big((\Delta \si_1)\nabla_a \si_2- (\Delta \si_2)\nabla_a \si_1+J(\si_2\nabla_a\si_1-\si_1\nabla_a \si_2)\big) \end{array}
  \right)
\end{equation} 
is a parallel tractor 2-form.
On the right $g$ is any metric from
the conformal class $[g_1]=[g_2]$ and 
 $\nd_b$ is the corresponding Levi-Civita connection.
By Proposition \ref{cKilltractor} the primary part
$$
k^b:=X_AZ_{B}{}^b\K^{AB}=\si_1\nd^b\si_2 -\si_2\nd^b \si_1
$$
 is a conformal Killing field. On the other hand since $\si_1$ is
nowhere vanishing we may work in the scale $g_1=\si_1^{-2}\bg$; with
$\nabla$ the corresponding Levi-Civita connection, we have
$\nabla\si_1=0$. Thus we see that $k^b$ is a gradient:
$$
k^b 
= g^{bc}_1\nabla_c s
$$
where $s$ is the function $\si_2/\si_1$. Note that away from the zero set 
of $s$, $g_2=s^{-2}g_1$ is Einstein.

$\Leftarrow :$ By definition $g_2$ is almost Einstein if and only if
$\nd_{(a}^{g_1}\nd^{g_1}_{b)_0}\si_2 + P^{g_1}_{(ab)_0}\si_2=0$ (see
\nn{prim}). Now recall that $\nd^{g_1}_a\si_1=0 $. Using this and that
$\si_1$ is non-vanishing, we see that $\nd_{(a}k_{b)_0}=0$ implies
$\nd_{(a}^{g_1}\nd^{g_1}_{b)_0}\si_2=0$.  On the other hand since
$g_1$ is Einstein $P^{g_1}_{(ab)_0}=0$. \quad $\Box$\\
\noindent {\bf Remark:} Note that not every conformal gradient comes
from a $\nabla$-parallel $2$-form tractor.  \endrk

By an obvious variation of the proof above we have the following.
\begin{proposition}\label{aeae}
If $(M,c,\si_1)$ is an almost Einstein structure then 
$(M,c,\si_2)$ is also an almost Einstein structure if and only
if the vector field
\begin{equation}\label{kdef2}
 k^a:=\si_1\nd^a\si_2 -\si_2\nd^a \si_1
\end{equation}
is a conformal Killing vector field. 
\end{proposition}
\noindent {\bf Remark:} Note that $k^a$ in \nn{kdef2}
is not a gradient if there is any point where both $\si_1$ and $\si_2$
vanish. One might term such a vector $k^a$ an almost conformal gradient.
\endrk

\subsection{Multiple almost Einstein structures: geometric implications}
\label{crae}
We will see from the explicit formula \nn{kpar} that the existence of
two linearly independent almost Einstein structures in a given fixed
conformal class has strong implications for the geometry and topology
of the underlying space. This uses some powerful results concerning
essential conformal vector fields.

\begin{theorem} \cite{AlekSb,Fe2,Yoshi} \label{Alekgl}
Let $(M^d,g)$, $d\geq 2$, be a Riemannian space, without boundary,
admitting a complete and essential conformal vector field. 
\begin{enumerate}
\item
If $(M,g)$ is non-compact then it is conformal to the 
Euclidean space $\RR^d$.
\item
If $(M,g)$ is compact then it is conformal to the sphere $S^d$ with
round metric $g_{rd}$.
\end{enumerate}
\end{theorem}
The case $d=2$ of the above is classically known from the theory of
Riemann surfaces.  Theorem \ref{Alekgl} says, in particular, that any
Riemannian space with a complete essential vector field is conformally
flat.

\begin{theorem} \label{zsp} Let $(M^d,g)$, $d\geq 2$, be a closed
  Riemannian manifold which admits a conformal gradient. Then
  $(M^d,g)$ is conformal to the round sphere $(S^d,g_{rd})$. 
\end{theorem}  
\noindent{\bf Proof:} Any gradient on $M$ has zeros, since any
function on a compact manifold has critical points. However by Theorem
\ref{kzeros} a conformal gradient field with a zero is
essential.  So the result follows from Theorem \ref{Alekgl}.
\quad $\Box$

\medskip

\begin{proposition}\label{oneRF}
  Suppose that $(M^{d\geq 3},c)$ admits two linearly independent
  parallel standard tractors $I_1$ and $I_2$, where $I_1$ is scalar
  flat and $\sigma_1=h(X,I_1)$ has a zero at $p\in M$.
Then  the primary part 
$k:= \sigma_1\nabla^a\sigma_2-\sigma_2\nabla^a\sigma_1$ of 
the parallel $\frac{1}{2}I_1\wedge I_2$ 
 is an essential conformal Killing field.\
\end{proposition}
\noindent{\bf Proof:}
Note that since $I_1$ and $I_2$ are linearly independent we have
$I_1\wedge I_2\neq 0$.  Since by Theorem \ref{classthm} (or
equivalently Theorem \ref{maina}) $j^1_p\si_1=0$ it follows from
\nn{kpar} that the $2$-form tractor $\frac{1}{2}I_1\wedge I_2$ is
given at $p\in M$ by
\[ \left( \begin{array}{ccc} & 0 &\\
\phantom{0}0&& \frac{1}{d} \si_2\Delta \si_1 \\
&\frac{1}{d}(\Delta \si_1)\nabla_a \si_2 & \end{array}
\right),
\] 
with respect to any metric $\tilde{g}\in c$.
This shows that $k$ has a
zero at $p$, and (by \nn{KKform})
either the divergence $\nabla_ck^c$ or $\nabla_{a}\nabla_ck^c $ 
does not vanish at $p$. 
Since $\tilde{g}$ was an arbitrary metric from the conformal class it
follows that the primary part $k$ of $\frac{1}{2}I_1\wedge I_2$ is an
essential conformal vector field. \quad $\Box$
 
Next note that if $I_1$ and $I_2$ are linearly independent scalar
negative almost Einstein structures such that zero sets of the scales
$\sigma_1$ and $\sigma_2$ intersect non-trivially then the projecting
part $k=\sigma_1\nabla^a\sigma_2-\sigma_2\nabla^a\sigma_1$ of
$\frac{1}{2}I_1\wedge I_2$ has a zero on $(M,c)$. However, the
vector $k$ is only an almost gradient field, in the sense of the
Remark in Section \ref{crE}, and so we cannot argue that $k$ is
essential.

Using our observations so far we led to the characterisation of
 compact Riemannian
spaces with at least $2$ linearly independent directed almost Einstein
structures as given in Theorem \ref{nonitsure}.\\ 
\noindent {\bf Proof of Theorem \ref{nonitsure}.}  Let $I_1$ and $I_2$ be two
linearly independent parallel tractors on $(M,c)$ and assume that
$I_1$ is Einstein that is $\si_1:=h(X,I_1)$ is nowhere vanishing.
Then, from Theorem \ref{2eincgrad}, the projecting part $k$ of
$\frac{1}{2}I_1\wedge I_2$ is a conformal gradient with respect to
$g_1:=\si_1^{-2}\bg$ and so the result follows from Theorem \ref{zsp}.

Now let us assume that $(M,c)$ is not conformally Einstein but that
there  exists a proper almost Ricci-flat structure, that is a
null tractor $I_1$ with $h(X,I_1)$ zero at some point. For any
linearly independent parallel standard tractor $I_2$ (from the class
assumed to exist on $(M,c)$) it follows from Proposition \ref{oneRF}
that the primary part of $\frac{1}{2}I_1\wedge I_2$ is essential.
From Theorem \ref{Alekgl} we can again conclude that $(M,c)$ is
conformal to the round sphere.

The remaining possibility is that every parallel standard tractor $K$
on $(M,c)$ is scalar negative and has a non-empty scale singularity
set $\S_K$. But then from Theorem \ref{classthm} (or equivalently
Theorem \ref{maina}) it follows that $\S_K$ is a totally umbillic
hypersurface.   \quad $\Box$

Note that a result by S.\ Gallot states that if the Riemannian
cone of a compact Riemannian space $(M,g)$ has decomposable holonomy
then $(M,g)$ is isometric to the round sphere \cite{Gallot}. 
 Observe that this is consistent with the
statement of Theorem \ref{nonitsure} as follows.  In the case of an
Einstein space with positive scalar curvature the Fefferman-Graham
ambient metric is just a product of the cone with a line
\cite{GrH1,LeitNorm}, and it follows (see \cite{GLsub}) that
almost Einstein structures are just parallel vectors on the cone. On
the other hand parallel vectors on the cone mean that the holonomy is
decomposable.
   
Before beginning the construction of examples we pause to prove the
other result from the introduction.\\
\noindent {\bf Proof of Theorem \ref{min}.} Since $K_1$ and $K_2$ are
parallel, and the tractor connection preserves the tractor metric, it
follows that $h(K_1,K_2)$ is constant on $M$ (at least if $M$ is
connected, as we assume for simplicity). On $M\setminus \Sigma_{K_1}$
$\si_1:=h(X,K_1)$ is non-vanishing and we may work in the metric scale
$\g=\si_1^{-2}\bg$. There and then we have 
$$
K_1\stackrel{\g}{=} (\si_1,~0,~-\frac{1}{d}\J^{\g}\si_1 ).
$$ On the other hand from \cite[Proposition 3.6]{Go-al} along
$\Sigma_{K_2}$, $K_2$ agrees with the normal tractor (of
\cite[Definition 2.8]{BEG})
$$
K_2|_{\Sigma_{K_2}}=N_{\Sigma_{K_2}} \stackrel{\g}{=} (0 ,~\nabla^{\g}\si_2,~-H^{\g} ),
$$ where $\si_2:=h(X,K_2)$ and $H^{\g}$ is the mean curvature of
$\Sigma_{K_2}$ in the metric scale $\g$. Now in fact $H^{\g}$ here is
the mean curvature as a weight $-1$ conformal density. The mean
curvature as a function in the scale $\g$ is $\underline{H}^{\g}=\si_1
H^{\g}$. But from the last two displays this is exactly (the
negative of) the tractor inner product of $K_1$ with
$N_{\Sigma_{K_2}}$. That is:
$$
\underline{H}^{\g}= -h( K_1, N_{\Sigma_{K_2}}).
$$ But, as remarked above, the conformally invariant inner product
$h(K_1,K_2)$ is constant on $M$, while evidently along $\Sigma_{K_2}$ we have 
$h(K_1,K_2)= h( K_1, N_{\Sigma_{K_2}})$. So we may conclude that 
$$
\underline{H}^{\g}=-h(K_1,K_2),
$$ and this is constant. Note that in particular, this constant is
zero if $K_1$ and $K_2$ are orthogonal.  \quad $\Box$\\
\noindent {\bf Remark:} Note that it is clear, from our considerations in the 
proof immediately above, that given an almost
Einstein structure $(M,c,K)$ and a smooth hypersurface $\Sigma$ with
normal tractor $N$, then $-h(K,N)$ gives a function along $\Sigma$
which generalises the mean curvature: it is defined on all of $\Sigma$
and is exactly the mean curvature $H^{\g}$ on $M\setminus
\Sigma_K$. Here $\Sigma_K$ is the zero set of $\si_K:=h(K,X)$ and
$\g=\si^{-2}_K \bg$. In fact this idea generalises beyond almost
Einstein spaces, but this will be taken up elsewhere.

Next observe that if we have $K_1$ and $K_2$ as in the Theorem
\ref{min} and we write $K_t:= \cos t K_1 +\sin t K_2$ then (away from
the zero set of $\si_1= h(X,K_1)$) for ${\mathbb R}\ni t \notin \pi{\mathbb Z}$ we have
$$
h(K_1,K_t)=-\underline{H}^{\g}_{\Sigma_{K_t}} 
$$ where $ \underline{H}^{\g}_{\Sigma_{K_t}}$ is the mean curvature of 
the zero set of $h(X,K_t)$ (at least for those values of $t$ such that $K_t$ is scalar negative) while for $t=0$
$$ 
h(K_1,K_0=K_1)=-S_{K_1}.
$$ So we see an interesting link between mean curvature and the scalar
curvature quantity $S_K$.  
\endrk

There are obvious results related to Theorem \ref{min}. For example in
the situation of that theorem $\S_{K_2}$ is itself directed almost
Einstein with parallel tractor given by the part of $K_1$ orthogonal
to $K_2$. This follows because from \cite[Theoerem 4.5]{Go-al}
the tractor connection on $(\S_{K_2},c_{K_2})$ is a restriction of the
ambient tractor connection to $K_2^\perp$. (Alternatively one may use
the conformal holonomy classification \cite{Armcomplete,LeitNCK} to
conclude that locally we have one of the special Einstein products
discussed below whence this also follows easily.) Thus in the setting
of Theorem \ref{min} we have the following, for example.
\begin{proposition}\label{minmore}
If $K_1$ and $K_2$ are spacelike and orthogonal then
$\S_{K_1}\cap \S_{K_2} $ is either empty or
else is a totally umbillic hypersurface in $(\S_{K_2},c_{K_2})$.
\end{proposition}
\noindent The reason we have taken $K_1$ and $K_2$ to be orthogonal is
that otherwise $\S_{K_1}\cap \S_{K_2} $ may be a (double) point,
collection thereof.

\subsection{Tractor bundles in low dimensions} \label{lowd} In
dimensions 1 and 2 a canonical tractor connection (equivalently Cartan
connection) is not determined locally by the conformal structure in
the usual way.  For some later constructions we use a notion of
tractor connection defined on (conformally) Einstein manifolds. In
dimension 2 a manifold is Einstein if and only if it has constant
scalar curvature. In dimension 1 any metric is Einstein.

Let $(M^d,c)$ be a conformal manifold. Observe that in dimensions
$d\geq 3$ if $\si^{-2}\bg=g\in c$ is Einstein, with corresponding
Einstein tractor $I$ (i.e.\ $I$ is parallel and $\si=h(X,I)$ is
nowhere vanishing), then the tractor connection is given on a section
of $\ce[1]\oplus \ce_a[1]\oplus \ce[-1]=[\cT]_g$ by
\begin{equation}\label{ndefE}
    \nd^\cT_a \left(\begin{array}{c} \alpha \\
  \tau^b \\ \rho \end{array}\right)
  \stackrel{g}{=} \left(\begin{array}{c} \nd_a \alpha - \tau_a \\
  \nd_a\tau^b + \delta_a{}^b( \rho +
  \mu \sigma^{-2}\alpha) \\ \nd_a\rho - \mu \si^{-2}\tau_a
  \end{array} \right) , \quad \mbox{or} \quad 
\left(\begin{array}{c} \nd_a s - t_a \\
  \nd_a t^b + \delta_a{}^b( r +
  \mu s) \\ \nd_a r - \mu t_a
  \end{array} \right)
\end{equation}
where $\mu=-\frac{1}{2}|I|^2= \underline{\Sc}^{g}/2d(d-1)$ and we
write $\underline{\Sc}^g$ to mean $g^{ab}\Ric^g_{ab}$ (so
$J^g/d=\si^{-2}\mu$). On the right we have simplified the connection
formula by trivialising the density bundles using the scale $\si$. So $s,r$
are simply functions and $t^a$ is an unweighted tangent field; these are
given by $s=\si^{-1}\alpha$, $\si \tau^a= t^a$, and $r=\si \rho$.
We will use this
formula to define the tractor connection in dimensions one and two as
follows. For simplicity we assume the manifold is connected. We
consider the case  of dimension two first. Let $(M,g)$ be a
Riemannian surface where $g$ is a constant scalar curvature metric.
In terms of the metric $g$ define $\ct\stackrel{g}{=}\ce[1]\oplus
\ce_a[1]\oplus \ce[-1]$ and equip this with the connection given by
the formula \nn{ndefE} above.  In dimension one the definition is the
same except that in this case we define such a connection for each
$\mu\in \mathbb{R}$ and $\mu$ is not related to the Levi-Civita
connection $\nabla$ or any intrinsic curvature quantity. The tractor
connections so defined for manifolds of dimension 2 are (locally)
flat.  The same is true in dimension 1 since in this case all
connections are locally flat.

\noindent{\bf Remarks:} It is straightforward to see that the
definition of the tractor connection above is equivalent to defining a
``Schouten tensor'' in dimensions 1 and 2.  For the Einstein metric
$g$ this is defined to be $P^g_{ab}:=\frac{1}{d}J^g\bg_{ab}$ where
$J^g:= \si^{-2}\mu d$.  Then it is defined for other metrics in the
conformal class by assuming the usual transformation formula
\nn{Rhotrans}.  With this definition one may then decree that the
tractor bundle also transforms according to the usual formula
\nn{transf} and use the usual formula \ref{ndef} for the tractor
connection (now using any metric from the conformal class).  That a
low dimensional conformal tractor bundle is equivalent to specifying a
Schouten tensor was pointed out in \cite{Cald}. The point we are
making here is that in dimension 2, when there is an Einstein metric
in the conformal class then this leads to a preferred choice.

Next we observe that there is a uniqueness result available for the
tractor connection of oriented (Riemannian signature) closed conformal
2-manifolds; in this dimension it draws on global data. By the
classical uniformisation result we know that there is a conformal
mapping to a metric of constant scalar curvature. On the other hand,
excepting the sphere, on such surfaces if there is a  constant scalar curvature
metric in the conformal class then it is unique see e.g.\
\cite{Mazzeo-Taylor} and references therein. Thus in such a conformal
class \nn{ndefE} gives the unique tractor connection.  \quad \endrk

Using this definition of the tractor connection it is straightforward
to see that the essential results from higher dimensions now extend to
dimensions $d=1$ and $d=2$. Evidently we have a preferred section
$X_A$ of $\ce_{A}[1]$ that gives a map $\ce^A \to \ce[1]$.  The
formula \nn{trmet} extends to give a metric $h$ in these dimensions by
working in the Einstein scale (or in other conformally related scales
according to the first Remark above).  This metric is then preserved
by the connection that \nn{ndefE} defines.  Clearly parallel sections
$I$ of the tractor bundle may again be divided into 3 classes
according to whether $S_I:=-h(I,I)$ is positive, zero, or negative. By
the same arguments as in higher dimensions (as in \cite{Go-al})
we then obtain the extension of Theorem \ref{classthm} to Einstein
manifolds of dimension 1 and 2. We shall exploit this in the
constructions below without further mention.

Although we shall not use it explicitly here, the tractor connection
in dimensions 1 and 2, as defined here, agrees with that induced from
the tractor connection on higher dimensional dimensional Einstein
manifolds by a suitable (local) totally umbillic embedding.  For
example for the unit round sphere $S^2$ we have $\mu=\frac{1}{2}$ and
it is easily verified that the connection $\nabla^\cT$ on $S^2$ given
by \nn{ndefE} agrees with the tractor connection induced on
$N^\perp\subset \cT|_\S$ along an equator $\S$ of the unit round
$S^3$. Here $N$ is the normal tractor of \cite[expression (12)]{BEG}
and writing $\cT_{S^3}$ for the tractor bundle of $S^3$, $N^\perp$ is
the subbundle of $\cT_{S^3}|_\S$ consisting of tractors orthogonal to
$N$.  Similarly, by an obvious extension of the definition of $N$ to
dimension 2, we then find that the $\mu=1/2$ connection on the unit
round $S^1$ agrees with the tractor connection induced on $N^\perp$
along an equator of $S^2$. It follows that the tractor connection on
$S^2$ and the $\mu=1/2$ tractor connection on $S^1$ have trivial
holonomy group. In fact for dimension 1 calculations it is convenient
to use the formula for the connection (also in the display
\ref{ndefE}) in terms of $s$, $t$ and $r$. It follows immediately from this 
that the equation of parallel transport in dimension 1 is equivalent to the ODE
\begin{equation}\label{ode}
(s''+ 2\mu s )'=0;
\end{equation}
parallel tractors are necessarily of the form $ (s,~t,~r) =
(s,~s',-(s''+\mu s))$ where $s$ solves the ODE.  So if $\mu=1/2$ then
we obtain $s=constant$ as well as a two dimensional family of
non-constant periodic solutions. So overall a three dimensional family
of solutions on $S^1$ as expected. On the other hand if we set
$\mu=-\frac{1}{2}$ then an $s=1$ solution shows that $(s,~t,~r) = (1,~
0,~ \frac{1}{2})$ is parallel, but since there are no periodic
non-trivial solutions to $(s''- s )=0$ it follows that that there are no global
parallel standard tractors on $S^1$ that are linearly independent of
$(1,~ 0,~ \frac{1}{2})$.

\subsection{Tractors on special Einstein products} 
The tractor bundle has rather special properties on certain product
manifolds. \\
\noindent{\bf Definition:} Let us say that $(M_1^{m_1}\times
M_2^{m_2},g_1\times g_2)$, with $m_1,m_2\geq 1$ and $m_1+m_2\geq 3$,
is a {\em special Einstein product} if
$(M_i,g_i)$, $i=1,2$, is Einstein with Einstein tractor $I_i$ and
$|I_1|^2=-|I_2|^2\neq 0 $.
By pull-back, we view each of
the bundles $\bT_i$, for $i=1,2$, as bundles on the product $M_1\times M_2$.
On any Einstein manifold, with
tractor bundle $\bT$ and parallel tractor $I$ giving the
Einstein structure, we shall write $\bT^\perp$ to mean the orthogonal
complement to $I$ in $\bT$.  Then we have (cf.
\cite{LeitNCK,Armcomplete}):
\begin{proposition}\label{prodT}
On a special Einstein product $(M_1\times M_2,g_1\times g_2)$ 
denote by $(\bT_i,h_i,\nabla^i)$, $i=1,2$, the respective 
standard tractor bundles.
Let us make the definitions: 
$\bT:=\bT_1^\perp\oplus \bT_2^\perp$ with bundle
metric $h=h_1\times h_2 $, where for $i=1,2$, $h_i$ is restricted to
$\bT_i^\perp$;
$\nabla:=\nabla^1\oplus \nabla^2$ where, for $i=1,2$, the $\nabla^i$
are also restricted to $\bT_i^\perp$. Then $(\bT,h,\nabla)$ is
the normal standard tractor bundle and connection.  The canonical
tractor $X\in \bT[1]$ is given by $X=(X_1^\perp , X_1^\perp
)\in \bT_1^\perp[1]\oplus \bT_2^\perp[1]$ where, for
$i=1,2$, $X_i^\perp$ is the orthogonal projection of $X_i$ onto $
\bT_i^\perp[1]$.
\end{proposition}
\noindent{\bf Proof:} Let $d:=m_1+m_2$.  Since the $I_i$, for $i=1,2$,
are not null, there are corresponding canonical orthogonal projections
$P^i:\bT_i\to \bT_i^\perp$ and each of these endomorphisms of $\bT_i$
is parallel for $\nabla^i$, $i=1,2$. We shall use the characterisation
of the normal tractor connection from \cite{CapGoTAMS}.

By construction $\bT$ is a vector bundle of dimension $d+2$ and the
metric $h$ is preserved by the restricted product connection
$\nabla$. If $I_1$ is spacelike ($|I_1|^2>0$) then $I_2$ is
timelike. Similarly if $I_1$ is timelike then $I_2$ is spacelike.  In
either case it is clear that $h$ has signature $(d+1,1)$. So the
connection $\nabla$ is a ${\frak g}$-connection where ${\frak
g}=so(d+1,1)$.

Trivialise the conformal density bundles via the metric $g:=g_1\times
g_2$ and view each of the canonical tractors $X_i$, $i=1,2$ as a
section of $\bT_i$.  Now define $X:= (X_1^\perp,X_2^\perp)$ where, for
$i=1,2$, $X_i^\perp$ is the image of $X_i$ under the orthogonal
projection to $\bT_i^\perp$, that is $X_i^\perp=P^i X_i$.  For $i=1,2$
we have that, in the scale $g_i$, $h^i(X_i,I_1)=1$ and thus from
$|I_1|^2=-|I_2|^2$ it follows that $X$ is null.  Working in the scale
$g_1\times g_2$, $X$ gives a map of $\ce$ into $\bT$ by $f\mapsto X
f$. Taking $f$ here as representing a section of $\ce[-1]$ we see
that, descending to the conformal class $[g_1\times g_2]$, $X$
determines a bundle map of $\ce[-1]$ into $\bT$ that we shall also
denote by $X$; the latter may be viewed as a section of $\ce[1]$. By
the construction of $\bT$, and using that $X$ is null we see that
$\bT$ has the composition structure $\bT=\ce[1]\lpl \ce_a[1]\lpl
\ce[-1]$. The projection $\bT\to \ce[1]$ is by $t\mapsto h(X,t)$ and
it is easily verified that restricting $h$ to the kernel
$\bT^{(0)}=\ce_a[1]\lpl \ce[-1] $ of this map recovers the conformal
metric $\bg$ for $(M,[g_1\times g_2])$ on $TM\cong
(\bT^{(0)}/\ce[-1])\otimes \ce[1]$.

Since, for $i=1,2$, $\nabla^i$ commutes with $P^i$, it follows that in
the scale $g_1\times g_2$ we have $\nabla X= (Z^1,Z^2)$ where the
$Z^i$ are the Z-projectors for each $[\bT_i]_{g_i}$ as in Section
\ref{trba}; we have used that in the given scale $P^i Z^i=Z^i$. It follows easily that
the connection $\nabla$ is non-degenerate in the sense of
\cite{CapGoTAMS} (this is the ``soldering'' condition).

Finally, since each of $I_1$ and $I_2$ is parallel, the curvature of
$\nabla$ is the direct sum of the tractor curvatures $\Om^1$ and
$\Om^2$ for, respectively, $\nabla^1$ and $\nabla^2$. Since
$\Om^iX^i=0$, i=1,2, it follows at once that it annihilates $X$. So
the curvature $\Om$ of $\nabla$ has no torsion and the leading
component of this is the Weyl curvature for $[g_1\times g_2]$. It
follows from the characterisation of \cite{CapGoTAMS} that the connection
$\nabla$ is the normal conformal standard tractor connection.  \quad
$\Box$\\

\noindent{\bf Remarks:} Note that an alternative proof of the above
follows easily via the cone product constructions of Fefferman-Graham
(ambient) metrics in \cite[Proposition 3.1]{GLsub} (and see also
Proposition 3.7 of that article) combined with the relationship
between the Fefferman-Graham metric and the conformal tractor calculus
as derived in \cite{CapGoamb}. One may also obtain the result from the
 formula \nn{ndef} by using that in the scale $g=g_1\times
g_2$ the Schouten tensor $P^g$ is the sum of the pull backs of the
Schouten tensors $P^{g_1}$ and $P^{g_2}$ from each component.

Note that in a special Einstein product the relationship between the
scalar curvatures (or rather lengths of the $I_i$) is exactly that
$(I_1,I_2)$ is null in $\bT_1\oplus \bT_2$ with respect to the product
tractor metric. \quad \endrk

\section{A (warped) product construction of almost Einstein manifolds}
\label{wprod}

\subsection{The general product construction} \label{gen}
It is immediate from Proposition \ref{prodT} that special Einstein
products are AE if and only if one of the factors has an additional AE
structure.
\begin{theorem}\label{prodAE}
Suppose that $(M,g)=(M_1\times M_2,g_1\times g_2)$ is a special
Einstein product. Then the conformal structure $(M,[g])$ admits a
(directed) almost Einstein structure if and if only one of the factors
$(M_1,[g_1])$ or $(M_2,[g_2])$ admits at least two linearly
independent (directed) almost Einstein structures.
\end{theorem}
\noindent{\bf Proof:} $I$ is a parallel non-zero standard tractor on
$(M,g)$ if and only if , for $i=1,2$, $P^i I $ are each parallel and
one of these is non-zero.  This proves the directed case. The general
result follows similarly and we leave this for the reader.  \quad
$\Box$\\ 
In fact the Theorem \ref{prodAE} above is the tractor
statement essentially equivalent to Proposition 3.4 from
\cite{GLsub}. Note also that in the theorem the generic case of
$(M,c)$ being a (directed) almost Einstein structure will be when both
of the factors admit at least two linearly independent (directed) almost
Einstein structures. Since in general a parallel tractor will project
into both parts ($\bT_1^\perp$ and $\bT_2^\perp$ in the notation of Proposition \ref{prodT}) of the tractor bundle.

Next we observe that the scale
singularity space is also easily identified.
\begin{corollary}\label{singprod}
Suppose that $(M,g)=(M_1\times M_2,g_1\times g_2)$ is a special
Einstein product.  Suppose that
$(M_1,g_1)$ has a
linearly independent set $\{K_1,\cdots, K_r\}$ of parallel tractors
which each take values in $\bT_1^\perp$. Then, via Proposition
\ref{prodT}, for each $i\in \{1,\cdots ,r\}$, $K_i$ determines the
almost Einstein structure $\tilde{K}_i:=(K_i,0)$ on $(M,g) $ and this
has the scale singularity set $\Sigma( \tilde{K}_i)= 
\Sigma_{K_i}\times M_2$, where $\Sigma_{K_i}\subset M_2$ is the scale
singularity set of $K_i$.
\end{corollary}
\noindent{\bf Proof:} This follows at once from Proposition
\ref{prodT} since $h(X,\tilde{K}_i)=\si_i$ where $\si_i$ is the pull back (via
obvious map $M=M_1\times M_2\to M_1$) of $h_1(X_1,K_i)$ to $M$.
\quad $\Box$\\
\noindent{\bf Remarks:} In relation to the last result, note in
particular that  $K_i$ is a proper almost Einstein structure (i.e.\
$\S_{K_i}\neq \emptyset$) if and only if  $\tilde{K}_i$ is proper.

Also observe that assuming the existence of a parallel standard
tractor $K$ which takes values in $\bT_1^\perp$ is the same as
assuming the existence of parallel tractor $K'$ so that $\{K',I_1\}$
is a linearly independent set; given such a tractor $K'$, $P^1 K'$
takes values in $\bT_1^\perp$.  \quad \endrk

Evidently the Corollary \ref{singprod} gives a general tool for
constructing proper almost Einstein structures. Let us set the
convention that in such products the first factor $(M_1,g_1)$ is
scalar positive, that is $S_{I_1}>0$. For the
construction of closed manifolds this way Theorem \ref{nonitsure}
severely limits the possibilities for $(M_1,g_1)$; to obtain a directed AE
structure $M_1$ is necessarily a sphere.  

\subsection{The examples $  S^{m_1\geq 2}\times  M^{m_2}_2$}\label{SM}
Let $(S^{m_1},g_1)$ be the standard unit sphere of dimension $m_1$. If
$m_1=2$ we fix $\mu=\frac{1}{2}$ in \nn{ndefE} to define the tractor
connection.  The sectional curvatures of this are 1, and so the
parallel tractor $I_1$ giving $g_1$ has $|I_1|^2=-1$. If
$(M_2^{m_2},g_2)$ is any (closed) Einstein manifold with $|I_2|^2=1$ (with $I_2$ the parallel tractor corresponding to $g_2$) then 
$\Ric^{g_2}=-(m_2-1)g_2$ and
$$
(M,g)=(S^{m_1}\times M_2,g_1\times g_2)
$$ is a (closed) special Einstein product. Via Corollary
\ref{singprod} the conformal manifold $(M,[g_1\times g_2])$ admits
$m_1+1$ linearly independent directed AE structures
$\{\tilde{K}_1,\cdots ,\tilde{K}_{m+1}\}$. By Theorem \ref{nonitsure}
these are necessarily proper and scalar negative. In fact the latter
point is already clear from Corollary \ref{singprod}: each non-zero
parallel tractor $K_i$ in $\mathcal{T}_1^\perp$ has $|K_i|^2>0$ since
$I_1$ is timelike. Thus the corresponding scale singularity set
$\S_{K_i}\subset S^{m_1}$ is a totally umbillic hypersphere. Then
since $|\tilde{K_i}|^2= |K_i|^2$, it is clear that $\tilde{K_i}$ is
scalar negative and its scale singularity set is the totally umbillic
hypersurface $\S_{\tilde{K}_i}=\S_{K_i}\times M_2$. By Theorem
\ref{min} each such hypersurface $\S_{\tilde{K}_i}$ is of constant
mean curvature with respect to the scale $\si^{-2}_{\tilde{K}_j}\bg$
determined by $\tilde{K}_j$ for $j\neq i$ in $\{1,\cdots ,m+1\}$; this
holds off $\S_{\tilde{K}_j}$ and the hypersurface is minimal if
$\tilde{K_i}$ and $\tilde{K_j}$ are orthogonal.

Note that if $K_1$ and $K_2$ are two linearly independent parallel
tractors in $I_1^\perp$ then it may be that the corresponding scale
singularity hyperspheres $\S_{K_1}$ and $\S_{K_2}$ intersect non-trivially. 
In fact generically $\S_{K_1}\cap\S_{K_2}\neq \emptyset$. On the other hand
$$ \S_{\tilde{K}_1}\cap\S_{\tilde{K}_2}=(\S_{K_1}\cap\S_{K_2})\times
M_1,
$$ and so generically this is also not empty. In any case where
$\S_{\tilde{K}_1}\cap\S_{\tilde{K}_2} \neq \emptyset$ the part
$M^{\geq 0}:=M\setminus M^-$ of $M$, where $\si_1:=h(X,\tilde{K_1})$
is nonnegative (the conformal density bundle $\ce[1]$ has a preferred
positive ray subbundle $\ce_+[1]$), is a Poincar\'{e}-Einstein
manifold. On this PE manifold $M^{\geq 0} $, with $\g:=\si_1^{-2}\bg$
on $M^+:= M^{\geq 0}\setminus \S_{K_1}$, the singularity hypersurface
of $K_2$ (viz.\ $\S_{K_2}\cap M^{\geq 0}$) is a totally umbillic
hypersurface which meets the conformal infinity of $(M^+,\g)$.  By
Proposition \ref{minmore} if $\tilde{K}_1$ and $\tilde{K}_2$ are
orthogonal then this intersection at is a totally umbillic
hypersurface of the conformal infinity.

Finally here we observe that one may be explicit about the almost
Einstein structures occurring in this class of cases. Let us consider
$ S^{m_1\geq 2}\times M^{m_2}_2$ to be embedded in $\mathbb{R}^{m_1+1}
\times M^{m_2}_2$ in the obvious way. Write $(x^1,\cdots ,x^{m_1+1})$
for the standard coordinates on the Euclidean space
$\mathbb{R}^{m_1+1}$.  Then away from $x^1=0$ on $ S^{m_1\geq 2}\times
M^{m_2}_2$ a typical metric in the class here is
$$
\g= (x^1)^{-2}(g_{\rm rd}\times g_2).
$$
Up to a conformal symmetry all the AE structures in this class are like this.

\subsection{The regular products $S^{1} \times M^{m_2}_2$ and periodic 
AE structures}\label{regs1} Since $S^1$ is 1-dimensional there is a
choice involved in the definition of the tractor connection. On the
standard unit circle let us term the connection \nn{ndefE} with
$\mu=\frac{1}{2}$ to be the {\em regular} tractor connection. Using
this the scale tractor $I_1$ on $S^1$ is parallel and has
$|I_1|^2=-1$. With reference to the ODE \nn{ode}, the parallel
tractors on $S^1$ orthogonal to $I_1$ are the solutions of $s''+s=0$
and so, taking the usual parametrisation, have two scale singularity
points.  Thus, as for the cases of $m_1\geq 2$ above, if
$(M_2^{m_2},g_2)$ is any (closed) Einstein manifold with
$\Ric^{g_2}=-(m_2-1)g_2$ then $ (M,g)=(S^1 \times M_2,g_1\times g_2) $
is a (closed) special Einstein product with a ($m_1+1=2$)-dimensional
space of almost Einstein structures compatible with the conformal
structure $[g_1\times g_2]$.  Everything works as in the case of
$m_1\geq 2$ above except that for two linearly independent parallel
tractors $K_1$ and $K_2$ on $S^1$ the corresponding scale singularity
sets, respectively, $\S_{K_1}$ and $\S_{K_2}$ do not intersect unless
they coincide. This has the obvious implication for the scale
singularity sets, $\S_{\tilde{K_1}}=\S_{K_1} \times M_1 $ and
$\S_{\tilde{K_2}}= \S_{K_2} \times M_1$, of the parallel tractors
$\tilde{K_1}$ and $ \tilde{K_2}$ on $(M,c)$.

Since the connections of \nn{ndefE} are natural for the given
Riemannian Einstein structure we may clearly pull these back to any
covering space. Thus for $K$ a parallel tractor on $S^1$ its pull back
(which we shall also denote by $K$) is a parallel tractor, and so a
directed almost Einstein structure, on any covering space of $S^1$.
It is straightforward to give a positive length parallel tractor $K$
on $S^1$ with two scale singularity points. Then on a $k$-fold
covering $\tilde{S}^1$ of $S^1$ the singularity set of $K$ consists of
$2k$ isolated points. The corresponding special Einstein product $
(M,g)=(M_1\times \tilde{S}^1,g_1\times g_2) $ (with $(M_1^{m_2},g_1)$
as above) is a periodic almost Einstein manifold with $2k$ scale
singularity points. This applies with $k$ infinite and $\tilde{S}^1$
is then $\mathbb{R}$. For finite $k$ the $2k$ singularity points may
be thought of as zeros of periodic solutions of (\ref{ode}) with appropriate 
$\mu$ in relation to the range of the independent variable.

\subsection{``Exotic'' $M^{m_1}_1\times S^{1}$ special Einstein products}
\label{exs1} 
If we equip the usual unit circle $(S^1,g_2)$ with the
tractor connection \nn{ndefE} with $\mu=-\frac{1}{2}$ then the
globally parallel tractor $I_2:
\stackrel{g_2}{=}(\si,~0,\frac{1}{2}\si^{-1}) $ has $|I_2|^2=1$, and so
in this sense $S^1$ behaves as a negative scalar curvature Einstein
manifold. (In this sense these cases are essentially a special case of those in Section \ref{SM} above. However some special comments are in order.) 
In particular we may form special Einstein products 
$$
(M,g):=(M^{m_1}_1\times S^{1},g_1\times g_2)
$$
where now 
$\Ric^{g_1}=(m_1-1)g_1$;  $g_1$ has positive scalar curvature.

There are no non-trivial parallel tractors on $S^1$ in $I_2^\perp$ and
so if $M_1$ is closed and Riemannian then it follows easily (using
also Theorem \ref{nonitsure}) that $(M,[g])$ admits compatible
directed AE structures if and only if $(M_1,g_1)$ is isometric to the
unit round sphere. Thus the only closed manifolds admitting directed
almost Einstein structures obtained in this way are the unit sphere
products as follows. (This shows we are really reduced to a special
case of the examples from Section \ref{SM}.)
\begin{proposition}\label{s1sn} The product of unit spheres 
$$
S^{n\geq 2}\times S^1
$$ with its usual metric and conformal structure admits an
$(n+1)$-dimensional vector space of parallel tractors, equivalently a
$(n+1)$-dimensional family of directed AE structures.  The product of
unit circles $S^1\times S^1$ admits a 4-dimensional family of directed
AE structures compatible with its given conformal structure.
\end{proposition}
\noindent Note that by the uniqueness (up to isomorphism) of the
normal tractor bundle it follows that the tractor bundle on the
product here is necessarily given by Theorem \ref{prodT}, with $S^1$
taken with the $\mu=-\frac{1}{2}$ tractor connection (as in
\nn{ndefE}).  The result for $S^1\times S^1$ follows as we may choose
which factor carries the positive length Einstein tractor.

\subsection{The examples $\mathbb{RP}_{m_1\geq 1} \times M^{m_2}_2$} 
\label{nond}
Recall that the standard conformal structure $c$ on
$\mathbb{RP}_{m_1\geq 1}$ admits an Einstein tractor $I_1$, with
$|I_1|^2=-1$ and in $I_1^\perp$ an $(m_1+1)$-dimensional family of
undirected almost Einstein structures.  We write $g_1$ for the
Einstein metric corresponding to $I_1$.  
We may
form products as for the cases above: If $(M_2^{m_2},g_2)$ is any
(closed) Einstein manifold with $\Ric^{g_2}=-(m_2-1)g_2$ then
$$
(M,g)=(\mathbb{RP}_{m_1}\times M_2,g_1\times g_2)
$$ is a (closed) special Einstein product.  Let us assume that $M_2$
is connected.  Since for each parallel tractor class $[K]$, with locally
$[K]\ni K' \in I_1^\perp$, and with singularity hypersurface
$\S_{[K]}$, the ``bulk'' $ \mathbb{RP}_{m_1} \setminus \S_{[K]}$ is
connected and so it follows that $M\setminus \S_{\tilde{[K]}}$ is also
connected, as recall $\S_{\tilde{[K]}} = \S_{[K]} \times M_2$ (here
locally $\tilde{K}=(K',0)$ as in Corollary \ref{singprod}). Thus
$\tilde{[K]} $ is undirected.

\section{Poincar\'e-Einstein subproducts and compactifications}\label{class}

\subsection{Decomposable conformal holonomy}
In this section we discuss scalar negative almost Einstein spaces
which have decomposable conformal holonomy.
This discussion will allow us to determine which of the 
Poincar{\'e}-Einstein collars, as introduced in \cite{GLsub}, are
suitable for compactification.

First we recall the notion of a metric cone.  Let $(M^d,g)$ be
a Riemannian space of dimension $d$. Then we define
$\hat{M}:=M\times\RR_+$. The cone metric of $g$ with timelike Euler
vector $t\cdot \partial_t$ is defined on $\hat{M}$ by
\[
\hat{g}\ :=\ -dt^2+ t^2\cdot g\ ,
\] 
where $t$ is the coordinate of $\RR_+$. We call $\hat{g}$ the timelike
cone metric of $g$. In \cite{GLsub} we have shown that any metric of
the form
\begin{equation}\label{dcneg}
g\ =\ r^{-2}(dr^2+(1-\mu r^2/2)^2 g_1+(1+\mu r^2/2)^2 g_2),
\end{equation}
where $\mu>0$ and $g_1,g_2$ are Einstein metrics of dimension $m_1\geq 1$,
resp., $m_2\geq 1$, with $2m_1(m_1-1)\mu:={\rm Sc}(g_1) $ and
$2m_2(m_2-1)\mu:=-{\rm Sc}(g_2) $, has the property that its timelike
cone $\hat{g}=-dt^2+t^2\cdot g$ is decomposable as a Riemannian
manifold.  
Note that with the coordinate change
$s=ln(\sqrt{\frac{\mu}{2}}\cdot r)$ for $r>0$ the metric $g$ takes the
form
\[
ds^2+2\mu( \sinh^2(s)\cdot g_1+\cosh^2(s)\cdot g_2 )\ .
\]

In the case that the Poincar\'e-Einstein collar metric \nn{dcneg} has
$g_1$ a round sphere metric then it is easily seen directly, or via
the examples of section \ref{SM}, that we may compactify the collar.
So our focus here is to determine whether there are other cases where
compactification is possible.

Spaces with a decomposable cone have been discussed earlier in, for
example, \cite{Gallot} and \cite{LeitNorm,LeitNCK}.  In general, one can show
that any metric with a decomposable timelike metric cone is (at least
on some open dense subset) locally isometric to a metric of the form
(\ref{dcneg}).  In fact, if the base metric is complete then the
following result is true.
\begin{theorem} \cite{devil} \label{ACG} Let $(M,g)$ be a complete
  Riemannian manifold of dimension at least two and such that its
  timelike cone $\hat{M}$, with cone metric $\hat{g}:=-dt^2+t^2g$, has
  decomposable holonomy algebra.  Then, on some open dense subset and for each 
connected component thereof,  $g$
  is  isometric to a metric of the form
\[
r^{-2}(dr^2+(1-\mu r^2/2)^2 g_1+(1+\mu r^2/2)^2 g_2), 
\]
with $\mu>0$, where the $g_2$ term is possibly absent,  and where
$g_1$ has positive constant sectional curvature. 
\end{theorem}

We shall also need the notion of conformal holonomy,
which we now recall. Let $(M^d,c)$ be a conformal
structure of dimension $d\geq 3$.  On $(M^d,c)$ we have the tractor
bundle $(\mathcal{T},\nabla)$ with its canonical normal connection. As for any
vector bundle connection, the tractor connection $\nabla$ on
$\mathcal{T}$ has a uniquely defined holonomy algebra and group (cf.\
\cite{KobNom}).  In the following, we denote the holonomy group of
$(\mathcal{T},\nabla)$ on $(M^d,c)$ by $Hol(\mathcal{T})$. We shall call
$Hol(\mathcal{T})$ the conformal holonomy group of $(M^d,c)$
(cf. \cite{Armcomplete,LeitNorm,LeitNCK}). The corresponding conformal
holonomy algebra is denoted by $\frak{hol}(\mathcal{T})$. The
conformal holonomy group is a subgroup of the M{\"o}bius group
$O(d+1,1)$, whose standard representation space is $\RR^{d+1,1}$.  If
the restriction of the standard representation of $Hol(\mathcal{T})$
on $\RR^{d+1,1}$ decomposes into a direct sum, then we call the
conformal holonomy group decomposable. Accordingly, if the conformal
holonomy algebra $\frak{hol}(\mathcal{T})$ splits into a direct sum,
then we call $\frak{hol}(\mathcal{T})$ decomposable.  Note that for
simply connected spaces $(M^d,c)$ decomposability of
$\frak{hol}(\mathcal{T})$ and $Hol(\mathcal{T})$ are equivalent.

  Furthermore, note that if the
conformal structure $c$ on $M^d$ includes an analytic metric $g$,
then the tractor connection $\nabla$ is analytic. This implies that
the conformal holonomy algebra $\frak{hol}(T)$ of $(M,c)$ is equal
to the conformal holonomy algebra of any open subspace $U$ of
$(M,c)$. Here we mean  $U$ is equipped with the restricted conformal
structure. It is also a  fact that the conformal holonomy
group and the holonomy group of the timelike cone over any Einstein
space of scalar curvature $-d(d-1)$ coincide (see \cite{GLsub}).  This
shows that any metric of the form (\ref{dcneg}) has decomposable
conformal holonomy group.

Now let us assume that $(M^d,c)$ with $dim(M)=d\geq 3$ is a compact
conformal space, which admits an almost Einstein structure $I$ which
is scalar negative (that is $S_I <0$),  and whose conformal holonomy
algebra $\frak{hol}(\mathcal{T})$ decomposes as a direct sum
$\frak{a}\times \frak{b}$ on $\RR^{d+1,1}$. Then we can describe the
local geometry of $c$ almost everywhere as follows.
\begin{theorem}\label{answer}
Let $(M^d,c,I)$ be either a Poincar\'e-Einstein manifold or a closed
scalar negative directed AE structure. Suppose further that $(M^d,c)$
has a decomposable conformal holonomy algebra
$\frak{hol}(\mathcal{T})$.  Then, on an open dense subset of $M$, the
conformal class $c$ locally includes a metric of the form
(\ref{dcneg}), where $g_1$ has positive constant sectional curvature.
\end{theorem}

\noindent{\sc Proof.} Let $I$ denote the given almost Einstein structure with
$S_I <0$; replacing $I$ by a constant multiple if necessary we may
assume $S_I=-1$.  We denote by $M^o$ the complement of the scale
singularity set $\Sigma_I$, of $I$ in $M$.  Writing $\si:=h(X,I)$,
then, in the restriction of $c$ to $M^o$, $g^o=\si^{-2}\bg$ is
complete Einstein metric of negative scalar curvature $-d(d-1)$. We
denote by $(\tilde{M},\tilde{g})$ the universal covering space of
$(M^o,g^o)$. Since, as a subset of $M$, the space $(M^o,g^o)$ has a
decomposable conformal holonomy algebra, the conformal holonomy group
of the universal covering space $(\tilde{M},\tilde{g})$ is
decomposable.  Hence we can conclude here that the timelike cone
$(\hat{M},\hat{g})$ over $(\tilde{M},\tilde{g})$ decomposes as a
Riemannian manifold. And since $(\tilde{M},\tilde{g})$ is complete, we
see via Theorem \ref{ACG} that the universal covering space takes
locally (on an open dense subset) the stated geometric form. But then
it is also clear that $(M^o,g^o)$ is locally isometric (on an open
dense subset) to a metric of the form (\ref{dcneg}), where $g_1$ has
positive constant sectional curvature.\quad  $\Box$\\

Theorem \ref{answer} has some immediate consequences.

\begin{corollary} \label{CorAA}
Let $(M^d,c)$ be a closed directed scalar negative AE space (resp.\ a
compact Poincar{\'e}-Einstein space) with decomposable conformal
holonomy algebra and with conformal infinity $(\Sigma,c_\S)$.  Then, on
an open dense subspace, $(\Sigma,c_\S)$ is locally conformally equivalent
to a 
product $S^{m_1}\times N$, where $S^{m_1}$ denotes the
unit sphere of dimension $m_1$, and $N$ is a space of dimension $m_2=d
-m_1-1$ equipped with an Einstein metric $h$ of negative scalar
curvature $-m_2(m_2-1)$.  Then, on an open dense subset of $M$, the space
$(M^d,c)$ is locally equivalent to the conformal structure of
$S^{m_1+1}\times N$ with metric $g_{rd}\times h$.
\end{corollary}
\noindent {\sc Proof.} Let us consider the conformal structure $c$ locally on a
 (small) collar of the conformal infinity  $(\Sigma,c_\S )$.  By
 Theorem \ref{answer} we know that the Poincar{\'e}-Einstein metric
 $g^o$, in the conformal class $c$, on such a collar is almost
 everywhere given by a metric of the form 
\[ r^{-2}(dr^2+(1-r^2/4)^2
 g_{rd}+(1+ r^2/4)^2 h)\ ,
\]
where $g_{rd}$ is the standard unit sphere metric (i.e.\ with constant
sectional curvature $1$ and $\mu=1/2)$. The displayed metric is
conformally equivalent to a metric $g_{rd}\times \tilde{h}$. Here the
coordinate $r$ is defined on some half interval $[0,\varepsilon)$.
This form of the Poincar{\'e}-Einstein metric also shows that the
(local) boundary conformal structure on $\Sigma$ is, on an open dense
subspace, given by $g_{rd}\times h$, where $g_{rd}$ is locally
isometric to the round unit sphere in dimension $m_1$. \hfill $\Box$

\begin{corollary} \label{CorBB}
Let $M^{m_1}\times N^{m_2}$ ($m_1\geq 2$, $m_2\geq 0$) with metric
$g=g_1\times g_2$ be either Einstein with $m_2=0$ or a special
Einstein product. Assume that, in either case, $(M,g_1)$ has positive scalar
curvature, but its sectional curvature is non-constant. Then the
Poincar{\'e}-Einstein collar of \cite{GLsub} for $g_1\times g_2$
cannot be compactified. That is, there is no conformal embedding of
that collar into a compact Poincar{\'e}-Einstein space.
\end{corollary}

\noindent {\sc Proof.}  Let us assume that there exists a 
 collar equipped with a Poincar{\'e}-Einstein metric of the form 
\[
r^{-2}(dr^2+(1-r^2/4)^2 g_1+(1+ r^2/4)^2 g_2)\ ,
\]
where $g_1$ has non-constant sectional curvature, which can be
embedded
into some compact Poincar{\'e}-Einstein space $(P,h)$ (whose conformal
infinity contains the boundary of that collar as a subspace). Now, the
conformal holonomy algebra of that collar is decomposable and it
follows, from the analyticity of the Einstein interior, that this
coincides with the conformal holonomy algebra of $(P,h)$.  
Thus, from  the
previous Corollary, on an open dense subspace,
the bulk of $(P,h)$ is locally isometric to a metric of the form
(\ref{dcneg}), with $g_1$ a piece of the round sphere. This 
contradicts the assumption on the embedded collar.  \hfill
$\Box$\\

Note that, in particular,  Corollary \ref{CorBB} states that the Poincar{\'e}-Einstein metric 
\[
\frac{1}{r^2}(\ dr^2\ +\ (1-\mu \cdot r^2/2)\cdot g_1\ )
\]
on some collar $[0,\varepsilon)\times M$ for an Einstein space
$(M^n,g_1)$ of dimension $n$ with positive scalar curvature cannot be
compactified, unless $(M,g_1)$ is conformally flat. In particular, if
$M$ is closed and simply connected then $(M,g_1)$ has to be (up to a
constant scale) the round unit sphere $S^n$.

\subsection{A global statement}\label{glob}

For each special Einstein product $g_1\times g_2$ Theorem 4.1 of
\cite{GLsub} gives a Poincar{\'e}-Einstein collar.  Using the
relationship between the conformal tractor connection and the ambient
metric, it is straightforward to see that by  construction, 
the Poincar{\'e}-Einstein collar, so obtained, 
 has the same conformal holonomy group as its boundary.  Both
holonomy groups are decomposable. In fact, it is also straightforward to
see that this equality of conformal holonomy groups of a boundary and
its bulk characterises, in the decomposable case, the
Poincar{\'e}-Einstein collar metric of \cite{GLsub} uniquely.

The following result states that Example 3.2 is the basic construction
principle for closed almost Einstein spaces, whose collar at the
singularity set is equivalent to the collar metric as introduced in
\cite{GLsub}.

\begin{theorem} \label{glob}
Assume that the space $(\Sigma,c_\S)$ is conformal to a special
Einstein product $A^r\times B^s$ of dimension $r+s=d-1$ with metric
$g_1\times g_2$, where $Sc(g_1)=r(r-1)$. If $(M^d,c,I)$ is a simply
connected, closed, and directed AE space with conformal infinity
$(\Sigma,c_\S)$ such that the conformal holonomy groups of
$(\Sigma,c_\S)$ and $(M^n,c)$ coincide, then $(A,g_1)$ is isometric to
$(S^r,g_{rd})$ and $(M^d,c)$ is conformal to $S^{r+1}\times B$ with
the metric $g_{rd}\times g_2$.
\end{theorem}
\noindent {\sc Proof.} The volume form $\alpha:=vol(g_2)$ for $g_2$ on
$B$ can be considered as a parallel $s$-form on $A\times B$ with
metric $g_1\times g_2$. This also means that $\alpha$ is a normalised
conformal Killing $s$-form on $(\Sigma,c_\S)$ in the sense of
\cite{LeitNCK,LeitNorm}, i.e., the corresponding $(s+1)$-form tractor
$K(\alpha)$ 
is
simple and parallel with respect to the tractor connection $\nabla$ on
$\Sigma$. The existence of this $(s+1)$-form tractor implies the
decomposability of the conformal holonomy group
$Hol(\Sigma,\mathcal{T})$. Since, by assumption,
$Hol(\Sigma,\mathcal{T})=Hol(M,\mathcal{T})$ and $M$ is simply
connected, it follows that there exists a parallel and simple
$(s+1)$-form tractor $T$ on $(M,c)$, whose restriction to
$\mathcal{T}$ on $\Sigma$ (which is in a natural and unique way a
subbundle of $\mathcal{T}$ on $M$, see \cite[Theorem 4.5]{Go-al})
gives rise to $K(\alpha)$.  We denote by $\beta$ the normalised
conformal Killing $s$-form on $(M^n,c)$, which corresponds to $T$.
Then, by construction, the pull back of the $s$-form $\beta$, by the
embedding $\iota:\Sigma\to M$, gives the $s$-form $\alpha$.

We want to argue now that the $s$-form $\beta$ has no zeros. This is
clear on the singularity set $\Sigma$, since $\iota^*\beta=\alpha$ has
no zeros.  So let us consider $\beta$ on the bulk $M^o$ of $M$ with
Poincar{\'e}-Einstein metric $\g$ and corresponding almost Einstein
structure $I$.  Recall that the product of a real line with the
timelike cone over $(M^o,\g)$ gives rise to the Fefferman-Graham
ambient metric of $(M^o,\g)$ \cite{LeitNCK,GLsub}.  Since, by
construction, $I\hook T=0$, it is clear that the $(s+1)$-form tractor
$T$ corresponds, via the Fefferman-Graham ambient metric, in a unique
way to a parallel and simple $(s+1)$-form $\hat{T}$ on the timelike
cone $\hat{M}$ over the bulk $(M^o,\g)$.  This $(s+1)$-form $\hat{T}$
has, by construction, negative length, i.e., it can be expressed by a
wedge product of one timelike $1$-form and some simple spacelike
$s$-form.  (Note that the $(s+1)$-form $\hat{T}$ can be interpreted as
the volume form of the timelike cone over $(B,g_2)$.)  
This shows that
inserting the Euler vector $t\cdot \partial_t$ into $\hat{T}$ is
nowhere vanishing. And, since the pull back of the $s$-form $(t\cdot
\partial_t)\hook \hat{T}$ via the embedding of $M^o$ (as the
$1$-level) into the timelike cone gives rise to $\beta$, we can
conclude that $\beta$ has no zeros on $M^o$.

Now we can choose a metric $h$ in $c$ such that $\beta$ has constant
length on $M$. In this scale $\beta$ is a parallel $s$-form
(cf. \cite{LeitNCK}).  Since $(M,h)$ is simply connected and closed, this
proves that $(M,c)$ is conformal to a special Einstein product
$C\times D$ with metric $h_1\times h_2$, where $Sc(h_1)>0$.  However,
from Corollary \ref{CorAA} we know that $g_1$ and $h_1$ have positive
constant sectional curvature. Moreover, since $C$ is closed and
simply connected it has to be the round sphere $S^{r+1}$.  And, since
$\hat{T}$ can be interpreted as the volume form of the timelike cone
over $(B,g_2)$, we conclude that $(D,h_2)$ is isometric to $(B,g_2)$.
Finally, since $I\hook T=0$, we see that $I$ on $M$ comes from an
almost Einstein structure on the factor $S^{r+1}$. This shows that the
boundary space $\Sigma$ is $S^r\times B$ with metric $g_1\times g_2$,
where $g_1$ is the standard sphere metric.  \hfill$\Box$\\

\end{document}